\input amstex
\documentstyle{amsppt}
\magnification=\magstep1
\hsize=5.2in
\vsize=6.8in
\topmatter
 \title COCYCLE SUPERRIGIDITY FOR PROFINITE ACTIONS\\
  OF PROPERTY (T) GROUPS \endtitle
 \rightheadtext{COCYCLE SUPERRIGIDITY}
\vskip 0.05in

\author {\rm ADRIAN IOANA}\endauthor

\address Math Dept  
Caltech, Pasadena, CA 91125\endaddress
\email aioana\@caltech.edu \endemail

\abstract
Consider a free ergodic measure preserving  profinite action $\Gamma\curvearrowright X$ (i.e. an inverse limit of actions $\Gamma\curvearrowright X_n$, with $X_n$ finite) of a countable property (T) group $\Gamma$ (more generally of a group $\Gamma$ which admits an infinite normal subgroup $\Gamma_0$ such that the inclusion $\Gamma_0\subset\Gamma$ has relative property (T) and $\Gamma/\Gamma_0$ is finitely generated)  on a standard probability space $X$. 
We prove that if $w:\Gamma\times X\rightarrow \Lambda$ is a measurable cocycle with values in a countable group $\Lambda$, then $w$ is cohomologous to a cocycle $w'$ which factors through the map $\Gamma\times X\rightarrow \Gamma\times X_n$, for some $n$.
As a corollary, we show that any orbit equivalence of $\Gamma\curvearrowright X$ with any free ergodic measure preserving action $\Lambda\curvearrowright Y$ comes from a (virtual) conjugacy of actions.
\endabstract
\endtopmatter
\document

\head \S {0. Introduction.}\endhead
During the past decade, the orbit equivalence theory of measure preserving actions of groups has been an extremely active area, with many new rigidity results having been proven (see the surveys [Sh1],[P4]). 
In particular,  certain classes of group--actions $\Gamma\curvearrowright X$ have been shown to be {\it orbit equivalent superrigid}, i.e. such that the equivalence relation $\Cal R_{\Gamma}$ on $X$ of belonging to the same $\Gamma$--orbit ($x\sim y$ iff $\Gamma x=\Gamma y$) remembers both the group $\Gamma$ and the action $\Gamma\curvearrowright X$ ([Fu1,2],[MSh1,2],[P1,2,3],[Ki1,2]). Since all ergodic actions  
$\Gamma\curvearrowright X$ of all infinite amenable groups $\Gamma$ induce isomorphic equivalence relations $\Cal R_{\Gamma}$ (up to a probability space isomorphism) ([Dy],[OW],[CFW]), such a rigidity phenomenon is very surprising and is characteristic to non-amenable groups only.

The main purpose of this paper is to present a new class of orbit equivalent superrigid actions. To explain this in more detail, we first review a few concepts, starting with the notion of orbit equivalence.  
Let $\Gamma\curvearrowright X$ be a free ergodic measure preserving action of a countable group $\Gamma$ on a standard probability space $(X,\mu)$ (i.e. isomorphic with the unit interval with the Lebesgue measure). Given another such action $\Lambda\curvearrowright Y$,  a probability space isomorphism $\theta:X\rightarrow Y$ is called an {\it orbit equivalence (OE)} between the actions $\Gamma\curvearrowright X$ and $\Lambda\curvearrowright Y$, if $\theta(\Gamma x)=\Lambda \theta(x)$, a.e. $x\in X$. If we moreover have that $\theta\Gamma{\theta}^{-1}=\Lambda$, then $\theta$ is called a {\it conjugacy} between $\Gamma\curvearrowright X$ and $\Lambda\curvearrowright Y$. 

Next, recall that an inclusion of countable groups $\Gamma_0\subset \Gamma$  has {\it relative property} (T) of Kazhdan--Margulis ([K],[Ma]) if any unitary representation of $\Gamma$ which weakly contains the trivial representation of $\Gamma$ must contain the trivial representation of $\Gamma_0$. Note that for $\Gamma=\Gamma_0$ this condition amounts to the {\it property} (T) of the group $\Gamma$.
Examples of relative property (T) inclusions of groups are given by $\Bbb Z^2\subset \Gamma\ltimes\Bbb Z^2$, for any non-amenable subgroup $\Gamma$ of SL$_2(\Bbb Z)$ ([Bu]) and by $\Gamma_0\subset \Gamma_0\times\Gamma_1$, for a property (T) group $\Gamma_0$ (e.g. $\Gamma_0=$SL$_n(\Bbb Z),n\geq 3$) and an arbitrary countable group $\Gamma_1$. 

Finally,  we call a measure preserving action $\Gamma\curvearrowright X$  {\it profinite} if it is the inverse limit of actions $\Gamma\curvearrowright X_n$, with $X_n$ finite probability spaces. For example, given a residually finite group $\Gamma$, let $G=\varprojlim \Gamma/\Gamma_n$ be the profinite completion of $\Gamma$ with  respect to a descending chain $\{\Gamma_n\}_n$ of finite index normal subgroups with trivial intersection. 
Then the (Haar) measure preserving action $\Gamma\curvearrowright G$ is free ergodic and profinite.

\proclaim {Theorem A (OE superrigidity)}  Let $\Gamma$ be a countable group with an infinite normal subgroup $\Gamma_0$ such that the inclusion $\Gamma_0\subset\Gamma$ has relative property (T) and $\Gamma/\Gamma_0$ is a finitely generated group. Assume that $\Gamma\curvearrowright^{\alpha}X$ is a free ergodic measure preserving profinite action on a standard probability space $X$ such that the restriction $\Gamma_0\curvearrowright^{\alpha_{|\Gamma_0}}X$ is also ergodic.
 
Let $\Lambda\curvearrowright^{\beta} Y$ be a free ergodic measure preserving action of a countable group $\Lambda$ on a standard probability space $Y$.  Suppose that $\theta:X\rightarrow Y$ is an orbit equivalence between $\alpha$ and $\beta$.  Then we can find an automorphism $\tau$ of $Y$ such that $\tau(y)\in\Lambda y$, a.e. $y\in Y$, two finite index subgroups $\Gamma_1\subset\Gamma$, $\Lambda_1\subset\Lambda$, a $\Gamma_1$--ergodic component $X_1\subset X$ and a $\Lambda_1$--ergodic component $Y_1\subset Y$ such that $(\tau\circ\theta)(X_1)=Y_1$ and $(\tau\circ\theta)_{|X_1}:X_1\rightarrow Y_1$ is a conjugacy between $\Gamma_1\curvearrowright^{\alpha_{|\Gamma_1}} X_1$ and $\Lambda_1\curvearrowright^{\beta_{|\Lambda_1}} Y_1$.
\endproclaim 

Before discussing some applications of Theorem A and its method of proof, let us give a brief history of previous results of this type.
The first such result was obtained by A. Furman, who combined Zimmer's cocycle superrigidity ([Z]) with ideas from geometric group theory to show that the actions SL$_n(\Bbb Z)\curvearrowright \Bbb T^n$ $(n\geq 3$) are OE superrigid ([Fu1,2]). 
This was followed by the work of N. Monod and Y. Shalom who employed techniques from bounded cohomology theory to prove that separately ergodic actions of products of hyperbolic groups are close to being OE superrigid ([MSh1,2]).
Recently, S. Popa used deformation/rigidity arguments in a von Neumann algebra framework to show that deformable actions (e.g. Bernoulli actions $\Gamma\curvearrowright [0,1]^{\Gamma}$) of rigid groups $\Gamma$ (e.g. groups $\Gamma$ which admit an infinite normal subgroups with relative property (T)) are OE superrigid ([P1,2], see also [PV]).
In subsequent work, Popa was able to remove the rigidity assumption on the group $\Gamma$ by assuming instead that $\Gamma$ is the product of two groups, one infinite and one non-amenable ([P3]).    
The last result along these lines is due to Y. Kida who proved that any ergodic action of any mapping class group is OE superrigid ([Ki1,2]).

As a consequence of Theorem A, we can construct uncountably many non--OE profinite actions for the arithmetic groups SL$_n(\Bbb Z)$ ($n\geq 3$), as well as for their finite index subgroups (see Corollary 5.3.), and for the groups SL$_m(\Bbb Z)\ltimes\Bbb Z^m$ ($m\geq 2$) (see Corollary 5.8.).  
In Section 4, we prove a more general version of Theorem A with stable orbit equivalence replacing orbit equivalence. This more general statement has the following interesting application: if $\Gamma\curvearrowright X$ is as above and if $X_0\subset X$ is a measurable set of irrational measure, then the equivalence relation $\Cal R_{\Gamma}^{X_0}:=\Cal R_{\Gamma}\cap (X_0\times X_0)$ cannot be induced by the free action of a countable group (compare with Theorem D in [Fu2] and Theorem 5.6. in [P2]). 

In proving Theorem A, we follow Zimmer's approach of studying orbit equivalences of actions via their associated OE cocycles. Recall in this respect that, given a measure preserving action $\Gamma\curvearrowright X$ and a countable group $\Lambda$, a measurable map $w:\Gamma\times X\rightarrow \Lambda$ is called a {\it cocycle} if it satisfies the relation $w(\gamma_1\gamma_2,x)=w(\gamma_1,\gamma_2 x)w(\gamma_2,x)$, for all $\gamma_1,\gamma_2\in\Gamma$ and a.e. $x\in X$. Two cocycles $w,w':\Gamma\times X\rightarrow \Lambda$ are {\it cohomologous} if there exists a measurable map $\phi:X\rightarrow\Lambda$ such that $w'(\gamma,x)=\phi(\gamma x)w(\gamma,x){\phi(x)}^{-1}$, for all $\gamma\in\Gamma$ and a.e. $x\in X$. Now, for an orbit equivalence $\theta:X\rightarrow Y$ between two free ergodic measure preserving actions $\Gamma\curvearrowright X$ and $\Lambda\curvearrowright Y$, consider the  measurable cocycle $w:\Gamma\times X\rightarrow \Lambda$ defined by the relation $\theta(\gamma x)=w(\gamma,x)\theta(x)$. A general principle then says that in order to show that $\theta$ is in fact a conjugacy, it essentially suffices to prove that $w$ is cohomologous to a homomorphism $\delta:\Gamma\rightarrow\Lambda$, i.e. to a cocycle which is independent of the $X$-variable ([Z],[Fu1,2],[P2]).
Guided by this principle we will deduce Theorem A as a consequence of the following:

 \proclaim {Theorem  B (Cocycle superrigidity)} Let $\Gamma\curvearrowright^{\alpha}X$ be as in Theorem A. Suppose that $\alpha$ is the limit of the actions $\Gamma\curvearrowright^{\alpha_n} X_n$, with $X_n$ finite, and, for every $n$, let $r_n:X\rightarrow X_n$ be the quotient map.
Let $\Lambda$ be a countable  group and $w:\Gamma\times X\rightarrow\Lambda$  be a measurable 
cocycle for $\alpha$.
Then there exists $n$ such that $w$ is cohomologous to a cocycle $w':\Gamma\times X\rightarrow \Lambda$ of the form $w'=w''\circ (id \times r_n)$, for some cocycle $w'':\Gamma\times X_n\rightarrow \Lambda$. 
\endproclaim 
In  other words: {\it any cocycle for $\alpha$ comes from one of the finite quotients $\alpha_n$}.	
On the other hand, notice that the $\alpha_n$'s do produce  non-trivial cocycles for $\alpha$.
Indeed, if we fix $n$ and identify $X_n$ (as a $\Gamma$-space) with $\Gamma/\Gamma_n$, for a finite index subgroup $\Gamma_n$  of $\Gamma$, then the natural cocycle $v_n:\Gamma\times\Gamma/\Gamma_n\rightarrow\Gamma_n$ lifts to a cocycle $v_n:\Gamma\times X\rightarrow \Gamma_n$.
 
The first examples of {\it cocycle superrigid} actions appeared only recently with the remarkable work of S. Popa who showed that if $\Gamma$ has an infinite normal subgroup with relative property (T), then any cocycle for the Bernoulli action of $\Gamma$  with values in any countable group $\Lambda$ (more generally, $\Lambda\in\Cal U_{fin}$)  is cohomologous to a homomorphism $\delta:\Gamma\rightarrow\Lambda$ ([P2])	.
Note that this result, in its most general form (see section 5 in [P2]), still requires that the action is weakly mixing. In contrast, Theorem B applies to profinite actions, which are, in some sense, farthest to being weakly mixing.   
It would be thus of interest to understand in a unitary way why cocycle superrigidity occurs for both Bernoulli and profinite actions.
More recently, Popa extended his result to cover many groups which do not have relative property (T) subgroups, such as $\Bbb F_2\times\Bbb F_2$ ([P3]).   
  
Related to the our results and motivated by the analogy with Popa's cocycle superrigidity results ([P2,3]), a natural open question arises: does the conclusion of Theorem B (or of Theorem A) hold true in the case $\Gamma=\Bbb F_2\times\Bbb F_2$? A positive answer to this question together with N. Ozawa and S. Popa's very recent work ([OP]) it would imply that any free ergodic measure preserving profinite action $\Gamma\curvearrowright X$ of $\Gamma=\Bbb F_2\times\Bbb F_2$ is {\it von Neumann superrigid}. This means that if $\Lambda\curvearrowright Y$ is another free ergodic measure preserving action such that the associated von Neumann factors ([MvN]) are isomorphic, $L^{\infty}X\rtimes\Gamma\simeq L^{\infty}Y\rtimes\Lambda$, then the actions $\Gamma\curvearrowright X$ and $\Lambda\curvearrowright Y$ are virtually conjugate, in the sense stated in Theorem A. 

The proof of Theorem B involves a ''zooming in and out'' argument, which we now briefly sketch. In the above context, assume for simplicity that $\Gamma$ has property (T) and let $w:\Gamma\times X\rightarrow \Lambda$ be a cocycle. Then, for a fixed $\gamma\in\Gamma$, as $n$ gets large,
we have that $(*)$ $r_n(x)=r_n(y)\Longrightarrow w(\gamma,x)=w(\gamma,y)$, with large probability. Next, using property (T) (for an appropriate representation of $\Gamma$), we deduce that $(*)$ holds true uniformly in $\gamma\in\Gamma$, as $n\rightarrow \infty$. Localizing, this tells us that for a large enough $n$, we can find $a\in X_n$ such that the map $X_a:=r_n^{-1}(\{a\})\ni x\rightarrow w(\gamma,x)\in\Lambda$ is almost constant, uniformly in $\gamma\in\Gamma_a$ (the stabilizer of $a$). Finally, a criterion for untwisting cocycles (see Section 2) implies that the restriction of $w$ to $\Gamma_a\times X_a$ is cohomologous to a homomorphism $\delta:\Gamma_a\rightarrow\Lambda$. Since $\Gamma_a$ has finite index in $\Gamma$, we can use this ''local'' information on $w$ to get the desired conclusion.  
  
In Section 1, we define the notion of profiniteness for actions and discuss certain properties of it.  Sections 3 and 4 are devoted to the proofs of theorems B and A, respectively, while in Section 5 we derive some applications of these results. In the last Section, we investigate two natural generalizations of profinite actions: compact  and weakly compact actions.

\vskip 0.1in
{\it Added in the proof.} After these results have been first circulated, A. Furman has been able to show
 that our main results still hold when the class of profinite actions is replaced by the larger class of compact actions (see Section 6 for more on compact actions).  

\head \S { 1. Profinite actions.}\endhead

\noindent
{\it Conventions.} All actions that we consider throughout this paper are of the form $\Gamma\curvearrowright^{\alpha} (X,\mu)$, where $\Gamma$ is a countable group which acts in a measure preserving way on a probability space $(X,\mu)$. For simplicity of notation, we will often denote an action by $\Gamma\curvearrowright X$. 
\vskip 0.05in
In this section we introduce the notion of profinite actions. After giving examples, we discuss ergodicity properties and freeness of such actions. We end the section by giving necessary and sufficient conditions for two profinite actions to be conjugate.
\vskip 0.05in
To define the notion of profinite action, we first recall the construction of an inverse limit of actions (see 6.3. in [G] for a reference). 
 Let $\Gamma\curvearrowright^{\alpha_n}(X_n,\mu_n)$ be a sequence of measure preserving  actions and assume that $\alpha_n$ is a quotient of $\alpha_{n+1}$, for all $n$. Let $q_n:(X_{n+1},\mu_{n+1})\rightarrow (X_n,\mu_n)$ be the quotient map, i.e. a measurable, measure preserving, onto map, such that $q_{n}(\gamma x)=\gamma q_n(x)$, a.e. $x\in X_{n+1}$, for all $\gamma\in\Gamma$. 
Define $$X=\{(x_n)_n|x_n\in X_n,q_{n}(x_{n+1})=x_n,\forall n\}$$ and let $r_n:X\rightarrow X_n$ be given by $r_n((x_m)_m)=x_n$. Then there exists a unique probability measure $\mu$ on $X$ such that $r_n:(X,\mu)\rightarrow (X_n,\mu_n)$ is measurable and measure preserving, for all $n$.
Denote by $\alpha$  the action of $\Gamma$ on $X$ given by $\gamma((x_n)_n)=(\gamma x_n)_n$. Then  $\alpha$ preserves $\mu$ and $r_n$ realizes $\alpha_n$ as a quotient of $\alpha$, for all $n$. We  say that  $\alpha$ is the {\bf limit} of $\alpha_n$ and we use the  notations $(X,\mu)=\varprojlim (X_n,\mu_n)$, $\alpha=\varprojlim\alpha_n$. Note that $\alpha$ is ergodic {\it iff} $\alpha_n$ is ergodic, for all $n$.
\vskip 0.1in
\noindent
{\bf 1.1. Definition.} A measure preserving action $\Gamma\curvearrowright^{\alpha}(X,\mu)$ is called {\bf profinite} if $\alpha=\varprojlim\alpha_n$, for a sequence of measure preserving actions $\Gamma\curvearrowright^{\alpha_n}(X_n,\mu_n)$ with $X_n$ finite.
\vskip 0.1in
\noindent
{\bf 1.2. Example.} Let $\Gamma$ be a countable group and let $\{\Gamma_n\}_n$ be a descending chain of finite index subgroups of $\Gamma$. For every $n$,  endow the quotient  $\Gamma/\Gamma_n$ with the  counting probability measure $\mu_n$ and denote by $\alpha_n$ the (transitive) left action of $\Gamma$ on the right cosets $\Gamma/\Gamma_n$. Also, for every $n$, let $q_n:\Gamma/\Gamma_{n+1}\rightarrow\Gamma/\Gamma_n$ be the map given by inclusion of cosets,  i.e. $q_n(x\Gamma_{n+1})=y\Gamma_n$ {\it iff} $x\Gamma_{n+1}\subset y\Gamma_n.$ Then $q_n$ is measure preserving and makes $\alpha_{n}$ a quotient of $\alpha_{n+1}$, for all $n$.
The limit action $\Gamma\curvearrowright^{\alpha}\varprojlim (\Gamma/\Gamma_n,\mu_n)$  is ergodic and profinite. 
In the case $\Gamma_n$ is a normal subgroup of $\Gamma$, for all $n\geq 0$, $\alpha$ can be alternatively seen as the action $\Gamma\curvearrowright (G,\mu)$, where $G$ is the profinite completion of $\Gamma$ with respect to $\{\Gamma_n\}_n$ and $\mu$ is the Haar measure of $G$.
\vskip 0.1in
\noindent
{\bf 1.3. Remarks.}
(1). Any ergodic, profinite action $\Gamma\curvearrowright^{\alpha} (X,\mu)=\varprojlim (X_n,\mu_n)$ arises as in the above example.  To see this, let $q_n:X_{n+1}\rightarrow X_n$ be the quotient maps and
 fix a sequence $a_n\in X_n$ such that $q_n(a_{n+1})=a_n$, for all $n$. Since $\alpha$ is ergodic, $\Gamma$ acts transitively on $X_n$, thus $X_n=\Gamma a_n$, for all $n$. Let $\Gamma_n=\{\gamma\in\Gamma|\gamma a_n=a_n\}$, then $\Gamma_{n+1}\subset\Gamma_n$. It is then easy to see that the map $\theta:X\rightarrow\varprojlim \Gamma/\Gamma_n$ given by $\theta((x_n)_n)=(\gamma_n\Gamma_n)_n$, where $\gamma_n$ is defined by the relation $x_n=\gamma_n a_n$, for all $x=(x_n)_n\in X$, is a probability space isomorphism identifying the actions $\Gamma\curvearrowright X$ and $\Gamma\curvearrowright\varprojlim\Gamma/\Gamma_n$.
\vskip 0.03in

(2). 
The representation of a profinite action as an action of the form $\Gamma\curvearrowright \varprojlim\Gamma/\Gamma_n$ is not unique.
Indeed, if $g_n$ is a sequence of elements of $\Gamma$ such that $g_n\Gamma_n\supset g_{n+1}\Gamma_{n+1}$, for all $n$, then the actions $\Gamma\curvearrowright \varprojlim\Gamma/\Gamma_n$ and $\Gamma\curvearrowright \varprojlim\Gamma/g_n\Gamma_ng_n^{-1}$ are isomorphic. This fact follows by applying the previous remark to $a_n=g_n\Gamma_n\in X_n=\Gamma/\Gamma_n$.

Also, we note that if $\{\Gamma_n^1\}_n$ and $\{\Gamma_n^2\}_n$ are two descending chains of finite index subgroups of $\Gamma$ such that $\Gamma_{n_{k+1}}^1\subset\Gamma_{N_k}^2\subset \Gamma_{n_k}^1$, for two subsequences $n_k$ and $N_k$, then the profinite actions $\Gamma\curvearrowright \varprojlim\Gamma/\Gamma_n^1$ and $\Gamma\curvearrowright\varprojlim\Gamma/\Gamma_n^2$ are isomorphic. 
 
\vskip 0.03in
(3).
Assume that $\Gamma\curvearrowright^{\alpha}(X,\mu)$ is an ergodic, profinite action.  
For every $m$, let $r_m:X\rightarrow X_m$ be the quotient map and denote $X_{a,m}=r_m^{-1}(\{a\})$, $\mu_{a,m}=\mu(X_{a,m})^{-1}\mu_{|X_{a,m}}$, for all $a\in X_m$. 
Also,  set $\Gamma_{a,m}=\{\gamma\in\Gamma|\gamma X_{a,m}=X_{a,m}\}=\{\gamma\in\Gamma|\gamma a=a\}$.
 Then the action $\Gamma_{a,m}\curvearrowright(X_{a,m},\mu_{a,m})$ is ergodic and profinite, for every $m$ and for all $a\in X_m$. 

 By the first remark we can assume that $\alpha$ is of the form $\Gamma\curvearrowright \varprojlim\Gamma/\Gamma_n$, for some descending chain $\{\Gamma_n\}_n\subset\Gamma$ of finite index subgroups,  and that $a=\Gamma_me\in\Gamma/\Gamma_m$. Under  this identification, the action  
$\Gamma_{a,m}\curvearrowright X_{a,m}$ becomes the action $\Gamma_m\curvearrowright\varprojlim_{n\geq m}\Gamma_m/\Gamma_n$, thus it is ergodic and profinite. 
\vskip 0.1in

For the next two results and their proofs we assume the  notations of Remark 1.3.(3). 

\proclaim{1.4. Lemma} Let $\Gamma\curvearrowright^{\alpha}(X,\mu)$ be an ergodic profinite action. Let $A\subset X$ be a measurable set which is  $\Gamma_0$-invariant, for some finite index subgroup $\Gamma_0$ of $\Gamma$. Then we can find $n$ and a $\Gamma_0$-invariant set $F\subset X_n$ such that $A=\cup_{a\in F}X_{a,n}.$
\endproclaim

{\it Proof.}
Assume first that $\Gamma_0$ is moreover a normal subgroup of $\Gamma$ and 
let $A\subset X$ be a $\Gamma_0-$invariant set. Let $B\subset X$ be an ergodic component for $\alpha_{|\Gamma_0}$. Since $\Gamma_0\subset\Gamma$ is normal we get that $\gamma B$ is an ergodic component, for all $\gamma\in\Gamma.$ Also, since $\alpha$ is ergodic, we have that $\cup_{\gamma\in S}(\gamma B)=X$, where $S\subset \Gamma$ satisfies $\Gamma=\sqcup_{\gamma\in S}(\gamma \Gamma_0)$. 
Since $A$ is $\Gamma_0$-invariant, by the above   facts we get that  $A=\cup_{\gamma\in S'}\gamma B$, for some subset $S'$ of $S$. Thus, it is sufficient to prove the conclusion in the case $A$ is an ergodic component for  $\alpha_{|\Gamma_0}$.  Indeed, if $B=\cup_{a\in F}X_{a,n}$, for some $n$ and $F\subset X_n$, then $A=\cup_{a\in S'F}X_{a,n}$.

\vskip 0.02in
Let $A$ be an ergodic component for $\alpha_{|\Gamma_0}$, then  $$\mu(\gamma A\cap A)\in\{0,\mu(A)\},\forall\gamma\in\Gamma\tag 1.4.a.$$

 From the way $(X,\mu)$ is constructed, we have that $\mu(A\Delta(\cup_{a\in F}X_{a,n}))<\mu(A)$, for some  $n$ and $F\subset X_n$. 
This implies that $\mu(A\cap(\cup_{a\in F}X_{a,n}))>\mu(\cup_{a\in F}X_{a,n})/2$, thus we can find $a\in F$ such that $\mu(A\cap X_{a,n})>\mu(X_{a,n})/2.$

 We claim that $A\cap X_{a,n}$ is $\Gamma_{a,n}-$invariant. Let $\gamma\in\Gamma_{a,n}$, then using (1.4.a.) and the fact that $\gamma X_{a,n}=X_{a,n}$, we get that $$\mu((A\cap X_{a,n})\cap \gamma(A\cap X_{a,n}))=\mu((A\cap\gamma A)\cap X_{a,n})\in\{0,\mu(A\cap X_{a,n})\}\tag 1.4.b.$$ On the other hand, we have that $$\mu((A\cap X_{a,n})\cap \gamma(A\cap X_{a,n,}))\geq 2\mu(A\cap X_{a,n})-\mu(X_{a,n})>0,\forall\gamma\in\Gamma_{a,n}\tag 1.4.c.$$  By combining (1.4.b.) and (1.4.c.) we deduce that $A\cap X_{a,n}$ is indeed $\Gamma_{a,n}-$invariant and
since  the action $\Gamma_{a,n}\curvearrowright X_{a,n}$ is ergodic (by Remark 1.3.(3)) we deduce that $A\cap X_{a\,n}=X_{a,n}$, hence that $X_{a,n}\subset A$, a.e. 
\vskip 0.02in
Next, we  prove that if $b\in X_n$ satisfies $\mu(A\cap X_{b,n})>0$, then $X_{b,n}\subset A$. For such   $b\in X_n$, let $\gamma\in\Gamma$ such that $b=\gamma a$. Then $X_{b,n}=\gamma X_{a,n}$ and  we have that $$0<\mu(A\cap X_{b,n})=\mu(\gamma^{-1}A\cap X_{a,n})\leq\mu(\gamma^{-1}A\cap A).$$  By using (1.4.a.)  this further implies that $\gamma^{-1}A=A$, thus $X_{b,n}=\gamma X_{a,n}\subset\gamma A=A$ a.e. 

In general, if we replace $\Gamma_0$ with the normal subgroup $\cap_{\gamma\in\Gamma}\gamma\Gamma_0\gamma^{-1}$ of $\Gamma$, then the above proof shows that $X=\cup_{a\in X}X_{a,n}$ for some $n$ and some set $F\subset X_n$. Since $A$ is $\Gamma_0-$invariant,  $F$ must also be $\Gamma_0-$invariant.
\hfill$\blacksquare$

\vskip 0.1in
\proclaim{1.5. Corollary} Let $\Gamma\curvearrowright^{\alpha} (X,\mu)$ be an ergodic, profinite action. Let $\Lambda$ be a countable group on which $\Gamma$ acts (e.g. assume that the action $\Gamma\curvearrowright \Lambda$ is given by $\gamma\cdot\lambda=\theta_1(\gamma)\lambda\theta_2(\gamma)^{-1}$, for two group homomorphisms $\theta_1,\theta_2:\Gamma\rightarrow\Lambda$).

  Let $\phi:X\rightarrow \Lambda$ is a measurable map such that $\phi(\gamma x)=\gamma\cdot\phi(x),$ for all $\gamma\in\Gamma$ and a.e. $x\in X.$ Then there exists   $n$ such that $\phi_{|X_{a,n}}$ is constant, for all $a\in X_n$.
\endproclaim
{\it Proof.} Let $\lambda\in\Lambda$, such that $X_{\lambda}=\{x\in X|\phi(x)=\lambda\}$ has $\mu(X_{\lambda})>0$. From the hypothesis we have that $\gamma X_{\lambda}=X_{\gamma\cdot\lambda},$ for all $\gamma\in\Gamma.$ This  implies  that $\mu(\gamma X_{\lambda}\cap 
X_{\lambda})\in\{0,\mu(X_{\lambda})\},$ for all $\gamma\in\Gamma$, and by applying the Lemma 1.4., we deduce that there exists $n\geq 0$ and $F\subset X_n$ such that $X_{\lambda}=\cup_{a\in F}X_{a,n}$.

Moreover, the ergodicity of $\alpha$ implies that $\cup_{\gamma\in\Gamma}(\gamma X_{\lambda})=X$, hence $\cup_{\gamma\in\Gamma}X_{\gamma\cdot\lambda}=X.$ Thus, if $\lambda'\in\Lambda$ satisfies $\mu(X_{\lambda'})>0$, then  $\lambda'=\gamma\cdot\lambda$, for some $\gamma\in\Gamma$, hence $$X_{\lambda'}=X_{\gamma\cdot\lambda}=\gamma X_{\lambda}=\cup_{b\in\gamma F}X_{b,n}.$$

\vskip 0.1in

\noindent
{\bf 1.6. Freeness of profinite actions.} Let $\Gamma\curvearrowright (X,\mu)=\varprojlim (\Gamma\curvearrowright(X_n,\mu_n))$ be a measure preserving ergodic profinite action with the quotient maps $r_n:X\rightarrow X_n$, for every $n$. For $\gamma\in\Gamma$ we have that $\{x\in X|\gamma x=x\}=\cap_n r_n^{-1}(\{x\in X_n|\gamma x=x\})$, hence $$\mu(\{x\in X|\gamma x=x\})=\lim_{n\rightarrow\infty} \mu (r_n^{-1}(\{x\in X_n|\gamma x=x\}))=$$ $$\lim_{n\rightarrow\infty}\mu_n(\{x\in X_n|\gamma x=x\})=\lim_{n\rightarrow\infty}|\{x\in X_n|\gamma x=x\}|/|X_n|.$$ 

Thus, the  action $\Gamma\curvearrowright (X,\mu)$ is {\bf (essentially) free} {\it iff} $$\lim_{n\rightarrow \infty}|\{x\in X_n|\gamma x=x\}|/|X_n|=0,\forall\gamma\in\Gamma\setminus\{e\}\tag 1.6.$$ 

In particular, if $\Gamma\curvearrowright (X,\mu)$ is free and if $\Gamma_n:=\{\gamma\in\Gamma|\gamma x=x,\forall x\in X_n\}$, then $\Gamma_n$ are finite index subgroups of $\Gamma$ and $\cap_n\Gamma_n=\{e\}$. Thus, any group $\Gamma$ which admits a free ergodic profinite action must be residually finite, a fact which we assumed implicitely in the statements of Theorems A and B.
Conversely, if $\Gamma$ is a residually finite group and if
$\{\Gamma_n\}_n$ is a descending chain of finite index, normal subgroups of $\Gamma$ with  $\cap_n\Gamma_n=\{e\}$, then by using (1.6.) it readily follows that the ergodic profinite action $\Gamma\curvearrowright \varprojlim(\Gamma/\Gamma_n,\mu_n)$ is also free.
\vskip 0.05in
Next, we consider a second construction of profinite actions where (1.6.) and thus freeness can be easily checked.
Let $\Gamma$ be a countable group which decomposes as a semidirect product $\Gamma=\Delta\ltimes\Gamma_0$. Assume $\{\Gamma_0^n\}_n$ is a descending chain of finite index subgroups of $\Gamma_0$ such that
$\Gamma_0^n$ is normal in $\Gamma$, for all $n$.
For every $n$ and all $a\in \Delta$, let $\theta_a^n\in$ Aut($\Gamma_0/\Gamma_0^n)$ be given by $\theta_a^n(c\Gamma_0^n)=(aca^{-1})\Gamma_0^n$, for all $c\Gamma_0^n\in\Gamma_0/\Gamma_0^n$.
Then $\Gamma$ acts on $\Gamma_0/\Gamma_0^n$ by the formula $$(a,b)\circ(c\Gamma_0^n)=\theta_a^n(bc\Gamma_0^n)=(abca^{-1})\Gamma_0^n,$$ for all $a\in\Delta,b\in\Gamma_0$, $ c\Gamma_0^n\in \Gamma_0/\Gamma_0^n$ and this action preserves the counting probability measure $\mu_n$ on $\Gamma_0/\Gamma_0^n$. If $q_n:\Gamma_0/\Gamma_0^{n+1}\rightarrow \Gamma_0/\Gamma_0^n$ is the map given by inclusion of right cosets, then $q_n$ is $\Gamma-$equivariant and measure preserving, for all $n$. Let $\Gamma\curvearrowright^{\alpha}\varprojlim (\Gamma_0/\Gamma_0^n,\mu_n)$ be the associated profinite action and note that, by construction,  $\alpha_{|\Gamma_0}$ is ergodic.

\proclaim {1.7. Lemma} In the above setting, assume 
that $\cap_n\Gamma_0^n=\{e\}$ and that, for all $a\in \Delta\setminus \{e\}$,
 the group $\{b\in\Gamma_0|[a,b]=e\}$ is of infinite index in $\Gamma_0$.
Then   $\alpha$ is  free.
\endproclaim

{\it Proof. } For every $n$ and $a\in\Delta$, denote Fix$(\theta_a^n)=\{x\in\Gamma_0/\Gamma_0^n|\theta_a^n(x)=x\}$. Then $q_n($Fix$(\theta_a^{n+1}))\subset$ Fix($\theta_a^n)$, thus $q_n$ induces a surjective homomorphism $$\overline{q}_n:(\Gamma_0/\Gamma_0^{n+1})/\text{Fix}(\theta_a^{n+1})\rightarrow (\Gamma_0/\Gamma_0^n)/\text{Fix}(\theta_a^n)\tag 1.7.a.$$

 We claim that $$\lim_{n\rightarrow\infty}[(\Gamma_0/\Gamma_0^n):\text{Fix}(\theta_a^n)]=\infty,\forall a\in\Delta\setminus\{e\}\tag 1.7.b.$$
If we assume  that this is not the case, then we can find $a\in\Delta\setminus\{e\}$ and $N$ such that $\overline{q}_n$ is an isomorphism, for all $n\geq N$.
 Let $\Gamma'=\{x\in\Gamma_0|\theta_a^N(x\Gamma_0^N)=x\Gamma_0^N\}$. Since $\overline{q}_n$ is injective, for all $n\geq N$, we get that $$\theta_a^n(x\Gamma_0^n)=x\Gamma_0^n,\forall x\in\Gamma',\forall n\geq N.$$ This rewrites $$x^{-1}axa^{-1}\in\Gamma_0^n,\forall x\in\Gamma',\forall n\geq N,$$ and since $\cap_{n}\Gamma_0^n=\{e\}$, we deduce that $[a,x]=e,$ for all $x\in\Gamma'$. However, since $\Gamma'\subset\Gamma_0$ is a finite index subgroup (as it contains $\Gamma_0^N$), this contradicts   the hypothesis, thus (1.7.b.) holds true.

Now, to check that $\alpha$ if free, let $(a,b)\in\Gamma\setminus\{e\}$, where $a\in\Delta$ and $b\in\Gamma_0$. Since the restriction $\alpha_{|\Gamma_0}$ is free, we can assume that $a\not=0$. We claim  that $$|\{x\Gamma_0^n\in\Gamma_0/\Gamma_0^n|(a,b)\circ (x\Gamma_0^n)=x\Gamma_0^n\}|\leq |\text{Fix}(\theta_a^n)|,\forall n\tag 1.7.c.$$ Indeed, this follows from the following fact: if $x,y$ verify $(a,b)\circ (x\Gamma_0^n)=x\Gamma_0^n$ and $(a,b)\circ (y\Gamma_0^n)=y\Gamma_0^n$, then $(y^{-1}x)\Gamma_0^n\in\text{Fix}(\theta_a^n)$.
By combining (1.7.b.) and (1.7.c.), we get that $$\lim_{n\rightarrow\infty}|\{x\Gamma_0^n\in\Gamma_0/\Gamma_0^n|(a,b)\circ (x\Gamma_0^n)=x\Gamma_0^n\}|/|\Gamma_0/\Gamma_0^n|=0,\forall a\in\Delta\setminus\{e\},\forall b\in\Gamma_0\tag 1.7.d.$$
Finally, by using  (1.6.) we deduce  that $\alpha$ is free.\hfill$\blacksquare$

\vskip 0.1in 
Recall that two measure preserving actions $\Gamma_i\curvearrowright^{\alpha_i}(X_i,\mu_i)$, $i=1,2$, are called {\bf conjugate} if there exists a group isomorphism $\delta:\Gamma_1\rightarrow\Gamma_2$ and a measure space isomorphism $\theta:(X_1,\mu_1)\rightarrow (X_2,\mu_2)$ such that $\theta(\gamma x)=\delta(\gamma)\theta(x)$, for all $\gamma\in\Gamma_1$ and a.e. $x\in X_1$. 
 
\proclaim {1.8. Proposition} Let $\Gamma_1$, $\Gamma_2$ be two countable groups and let $\{\Gamma_{1,n}\}_n\subset\Gamma_1$ and $\{\Gamma_{2,n}\}_n\subset\Gamma_2$ be descending chains of finite index subgroups.  Then the profinite actions $\Gamma_i\curvearrowright^{\alpha_i} (X_i,\mu_i):=\varprojlim (\Gamma_i/\Gamma_{i,n},\mu_{i,n})$, $i=1,2$, are conjugate {\it iff} there exists a group isomorphism $\delta:\Gamma_1\rightarrow\Gamma_2$, two subsequences $\{n_k\}_k,\{N_k\}_k\subset \Bbb N$ and $\gamma_{k}\in\Gamma_2$ such that 

 $$\gamma_{N_{k+1}}\Gamma_{2,N_{k+1}}\subset \gamma_{N_k}\Gamma_{2,N_k},$$ $$ \delta(\Gamma_{1,n_{k+1}})\subset\gamma_{N_k}\Gamma_{2,N_k}{\gamma_{N_k}}^{-1}\subset\delta(\Gamma_{1,n_k}),\forall k.$$

In particular, if for some $i\in\{1,2\}$ we have that $\Gamma_{i,n}$ is a normal subgroup of $\Gamma_i$, for all $n$,  then $\alpha_1$ is conjugate to $\alpha_2$ {\it iff} there exists a group isomorphism $\delta:\Gamma_1\rightarrow\Gamma_2$ and two subsequences $\{n_k\}_k,\{N_k\}_k\subset \Bbb N$ such that $$\delta(\Gamma_{1,n_{k+1}})\subset\Gamma_{2,N_k}\subset\delta(\Gamma_{1,n_k}),\forall k.$$

\endproclaim
{\it Proof.} 
Assume first that $\alpha_1$ is conjugate to $\alpha_2$. Let $\delta:\Gamma_1\rightarrow\Gamma_2$ be a group isomorphism and $\theta:(X_1,\mu_1)\rightarrow (X_2,\mu_2)$ be a measure space isomorphism implementing the conjugacy.
 For $i\in\{1,2\}$  and for all 
$n$, denote  $X_{i,n}=\Gamma_i/\Gamma_{i,n}$ and let $r_{i,n}:X_i\rightarrow X_{i,n}$ be the quotient map. Also, for all $a\in X_{i,n}$, denote $X_{i,a,n}=r_{i,n}^{-1}(\{a\})$ and $\Gamma_{i,a,n}=\{\gamma\in\Gamma
_i|\gamma X_{i,a,n}=X_{i,a,n}\}$.

\vskip 0.03in
{\it Claim.} For all $n$ and for every $a\in X_{1,n}$,  there exist $N$ and $F\subset X_{2,N}$ such that $\theta(X_{1,a,n})=\cup_{b\in F}X_{2,b,N}.$ Moreover, in this case, we have that $\Gamma_{2,b,N}\subset \delta(\Gamma_{1,a,n})$, for all $b\in F$.
\vskip 0.03in
{\it Proof of Claim}.
Since $X_{1,a,n}$ is $\Gamma_{1,a,n}$-invariant, we get that $\theta(X_{1,a,n})$  is $\delta(\Gamma_{1,a,n})$-invariant. 
 Using the fact that $\delta(\Gamma_{1,a,n})$ has finite index in $\Gamma_2$, Lemma 1.4.  implies that $\theta(X_{1,a,n})=\cup_{b\in F}X_{2,b,N},$ for some $N$ and some set $F\subset X_{2,N}$.
For the second assertion, let $b\in F$ and $\gamma\in \Gamma_2$ such that $\gamma X_{2,b,N}=X_{2,b,N}$. Thus, in particular, $\mu_2(\gamma\theta(X_{1,a,n})\cap \theta(X_{1,a,n}))>0$, hence $\mu_1(\delta^{-1}(\gamma)X_{1,a,n}\cap X_{1,a,n})>0$. Finally, this implies that $\delta^{-1}(\gamma)$ invaries $X_{1,a,n}$, which proves the claim.\hfill$\square$
\vskip 0.03in
  Similarly, we get that for all $n$ and every $b\in X_{2,n}$, there exist $N$ and $F\subset X_{1,N}$ such that $X_{2,b,n}=\cup_{a\in F}\theta(X_{1,a,N})$. Also, in this case, $\delta(\Gamma_{1,a,N})\subset \Gamma_{2,b,n}$, for all $a\in F$.
\vskip 0.03in
Now, for every $n$, define $a_n=e\Gamma_{1,n}\in X_{1,n}=\Gamma_1/\Gamma_{1,n}$. 
By the above claim it follows that for every $n$ we can find $m\geq n$, $N$ and $b\in X_{2,N}$ such that $\theta(X_{1,a_{m},m})\subset X_{2,b,N}\subset \theta(X_{1,a_n,n})$ and  $\delta(\Gamma_{1,a_{m},m})\subset \Gamma_{2,b,N}\subset\delta(\Gamma_{1,a_n,n}).$
Indeed, by applying the claim  we get that $\theta(X_{1,a_n,n})=\cup_{b\in F}X_{2,b,N}$ for some $N$ and some set $F\subset X_{2,N}$. Now, by applying the claim  a second time we can find $m\geq n$ such that for all $b\in F$, $X_{2,b,N}=\cup_{a\in F_b}\theta(X_{1,a,m})$, for some set $F_b\subset X_{1,m}$. Since we also have that $\theta(X_{1,a_m,m})\subset \theta(X_{1,a_n,n})$, we can find $b\in F$ such that $\theta(X_{1,a_m,m})\subset X_{2,b,N}$. 

Thus, we can inductively construct two subsequences $\{n_k\}_k,\{N_k\}_k\subset\Bbb N$ and $b_{k}\in X_{2,N_k}$ such that $$\theta(X_{1,a_{n_{k+1}},n_{k+1}})\subset X_{2,b_{k},N_k}\subset \theta(X_{1,a_{n_k},n_k})\tag 1.8.a.$$ and $$\delta(\Gamma_{1,a_{n_{k+1}},n_{k+1}})\subset \Gamma_{2,b_{k},N_k}\subset\delta(\Gamma_{1,a_{n_k},n_k}),\forall k\tag 1.8.b.$$

Now, for every $k$, let $\gamma_{k}\in\Gamma_2$ such that $X_{2,b_{k},N_k}=\gamma_{k}\Gamma_{2,N_k}$. By using inclusion (1.8.a.) we get that  $\gamma_{k+1}\Gamma_{2,N_{k+1}}\subset\gamma_{k}\Gamma_{2,N_k}$, for all $k$.
Since we also have that $\Gamma_{2,b_{k},N_k}=\gamma_{N_k}\Gamma_{2,N_k}{\gamma_{N_k}}^{-1}$ and that $\Gamma_{1,a_{n_k},n_k}=\Gamma_{1,n_k}$, for al $k$, inclusion (1.8.b.) gives the rest of the conclusion. 
\vskip 0.03in
Conversely, if we assume inclusions from the hypothesis to hold true, then Remark 1.3.(2) implies that $\alpha_1$ and $\alpha_2$ are conjugate.\hfill$\blacksquare$

\vskip 0.1in

\head \S {2. A criterion for untwisting cocycles.}\endhead

\vskip 0.1in

Let $\Gamma\curvearrowright^{\alpha}(X,\mu)$ be a measure preserving action and let $\Lambda$ be a Polish group.
 A measurable map  $w:\Gamma\times X\rightarrow \Lambda$ is called a {\bf cocycle} for $\alpha$   if it satisfies $$w(\gamma_1\gamma_2,x)=w(\gamma_1,\gamma_2x)w(\gamma_2,x),\forall \gamma_1,\gamma_2\in\Gamma$$ and a.e. $x\in X$. Two cocycles $w_1$ and $w_2$ are said to be {\bf cohomologous} if there exists a measurable map $\phi:X\rightarrow\Lambda$ such that $$w_1(\gamma,x)=\phi(\gamma x)w_2(\gamma,x)\phi(x)^{-1},\forall \gamma\in\Gamma$$ and a.e. $x\in X$. Note that if $\psi:\Gamma\rightarrow\Lambda$ is a group homomorphism, then $w(\gamma,x)=\psi(\gamma)$ defines a cocycle. 
In this Section, we give a criterion for untwisting cocycles, i.e. for showing that a cocycle is cohomologous to a homomorphism. This criterion will follow as a consequence of the next lemma whose proof is inspired by the proofs of 4.2. in [P2] and 4.2. in [Fu4].

\proclaim{2.1. Lemma} Let $\Gamma\curvearrowright^{\alpha}(X,\mu)$ be a  measure preserving ergodic action. Let $\Lambda$ be a countable  group and let $w_1,w_2:\Gamma\times X\rightarrow \Lambda$ be two cocycles for $\alpha$. Let $C\in (7/8,1)$ and assume that $\mu(\{x\in X|w_1(\gamma,x)=w_2(\gamma,x)\})\geq C,$ for all $\gamma\in\Gamma.$ Then $w_1$ is cohomologous to $w_2$.
\endproclaim
{\it Proof.} Endow $\Lambda$ with the counting measure $c$ and let $\sigma$ be the measure preserving action of $\Gamma$ on the infinite measure  space $(X\times\Lambda,\mu\times c)$ given by $$\gamma(x,\lambda)=(\gamma x,w_1(\gamma,x)\lambda{w_2(\gamma,x)}^{-1}),\forall\gamma\in\Gamma,\lambda\in\Lambda,x\in X.$$ Let $\pi:\Gamma\rightarrow\Cal U(L^2(X\times\Lambda,\mu\times c))$ be the induced unitary representation and let $\xi\in L^2(X\times\Lambda,\mu\times c)$ be given by $\xi(x,\lambda)=\delta_{\lambda,e}$. Then $$||\pi(\gamma)(\xi)-\xi||^2=2-2\Re <\pi(\gamma)(\xi),\xi>=$$  $$2-2\mu(\{x\in X|w_1(\gamma,x)=w_2(\gamma,x)\})\leq 2-2C,\forall\gamma\in\Gamma.$$ Thus $||\pi(\gamma)(\xi)-\xi||\leq \sqrt{2-2C},$ for all $\gamma\in\Gamma$, and a standard averaging argument shows that we can find $\eta\in L^2(X\times\Lambda,\mu\times c)$ such that $||\eta-\xi||\leq \sqrt{2-2C}<1/2$ and $\eta$ is $\pi(\Gamma)$-invariant, i.e. $$\eta(\gamma x,w_1(\gamma,x)\lambda w_2(\gamma,x)^{-1})=\eta(x,\lambda),\forall \gamma\in\Gamma,\lambda\in\Lambda, \tag 2.1.$$and a.e. $x\in X$.   

Let $X_0$ be the set of $x\in X$ such that there exists a unique $\lambda=\phi(x)\in\Lambda$ with $|\eta(x,\lambda)|>1/2$. 
Then identity $(2.1.)$ implies that $X_0$ is $\alpha(\Gamma)-$invariant and that $$w_1(\gamma,x)\phi(x){w_2(\gamma,x)}^{-1}=\phi(\gamma x),\forall\gamma\in\Gamma\tag 2.2.$$ a.e. $x\in X_0.$ 
On the other hand, the inequality $||\eta-\xi||<1/2$ is equivalent to $$\int_{X}[|\eta(x,e)-1|^2+\sum_{\gamma\in\Gamma\setminus\{e\}}|\eta(x,\gamma)|^2]d\mu(x)<1/4.$$
In particular, the set $$Y=\{x\in X||\eta(x,e)-1|^2+\sum_{\gamma\in\Gamma\setminus\{e\}}|\eta(x,\gamma)|^2<1/4\}$$ satisfies $\mu(Y)>0.$ It is then clear that $Y\subset X_0$, thus $\mu(X_0)>0$ and since $\alpha$ is ergodic, we get that $X_0=X$, a.e. Finally,  $(2.2.)$ implies that $w_1$ is cohomologous to $w_2$.   
\hfill$\blacksquare$

\proclaim {2.2. Corollary} Let $\Gamma\curvearrowright^{\alpha}(X,\mu)$ be a  measure preserving ergodic action. Let $\Lambda$ be a countable  group and let $w:\Gamma\times X\rightarrow \Lambda$ be a cocycle for $\alpha$. Let $C\in (7/8,1)$ and assume that $$(\mu\times\mu)(\{(x_1,x_2)\in X\times X|w(\gamma,x_1)=w(\gamma,x_2)\})\geq C,\forall\gamma\in\Gamma.$$ Then $w$ is cohomologous to a group homomorphism $\psi:\Gamma\rightarrow\Lambda$.
\endproclaim
{\it Proof.} For $\gamma\in\Gamma$ and $\lambda\in\Lambda$, denote $S_{\gamma,\lambda}=\{x\in X|w(\gamma,x)=\lambda\}$. If $\gamma\in\Gamma$, then since $\{S_{\gamma,\lambda}\}_{\lambda\in\Lambda}$ is a partition of $X$, we get that $\sum_{\lambda\in\Lambda}\mu(S_{\gamma,\lambda})=1$. Thus, $$C\leq (\mu\times\mu)(\{(x_1,x_2)\in X\times X|w(\gamma,x_1)=w(\gamma,x_2)\})=$$ $$\sum_{\lambda\in\Lambda}\mu(S_{\gamma,\lambda})^2\leq (\sum_{\lambda\in\Lambda}\mu(S_{\gamma,\lambda}))\max_{\lambda\in\Lambda}\mu(S_{\gamma,\lambda})=\max_{\lambda\in\Lambda}\mu(S_{\gamma,\lambda}),\forall\gamma\in\Gamma.$$ In conclusion, for every $\gamma\in\Gamma$, we can find a unique $\psi(\gamma)\in\Lambda$ such that the set $T_{\gamma}=\{x\in X|w(\gamma,x)=\psi(\gamma)\}$ has measure $\mu(T_{\gamma})\geq C$. Next, note that given $\gamma_1,\gamma_2\in\Gamma$, we have that

$$\psi(\gamma_1\gamma_2)=w(\gamma_1\gamma_2,x)=w(\gamma_1,\gamma_2x)w(\gamma_2,x)=\psi(\gamma_1)\psi(\gamma_2),\forall x\in T_{\gamma_1\gamma_2}\cap\gamma_2^{-1}T_{\gamma_1}\cap T_{\gamma_2}.$$

 Since $\mu(T_{\gamma_1\gamma_2}\cap\gamma_2^{-1}T_{\gamma_1}\cap T_{\gamma_2})\geq 3C-2>0,$ for all $\gamma_1,\gamma_2\in\Gamma$, we deduce that $\psi$ is a homomorphism and the claim follows by the previous lemma. \hfill$\blacksquare$

\head\S  {3. Proof of cocycle superrigidity.}\endhead
 In this Section we give the proof of Theorem B. To this end, recall that the action $\Gamma\curvearrowright^{\alpha} X$ is a limit of actions  $\Gamma\curvearrowright^{\alpha_n}X_n$, where $X_n$ are finite probability spaces together with $\Gamma$-equivariant quotient maps $r_n:X\rightarrow X_n$, for all $n\geq 0$.
For a fixed $n$, denote $X_{a,n}=r_n^{-1}(\{a\})$,  $\Gamma_{a,n}=\{\gamma\in\Gamma|\gamma a=a\}$, for all $a\in X_n$, and observe that $\mu(X_{a,n})=|X_n|^{-1}$. Also, recall that $\Gamma_0$ is a normal subgroup of $\Gamma$ such that the inclusion $\Gamma_0\subset\Gamma$ has relative property (T), the quotient $\Gamma/\Gamma_0$ is finitely generated and the restriction $\Gamma_0\curvearrowright^{\alpha_{|\Gamma_0}}X$ is ergodic.
\vskip 0.1in
Next, let $c$ be the counting measure on $\Lambda$ and  define the measure space $(Z,\rho)=(X\times X\times \Lambda,\mu\times\mu\times c)$. Let $\sigma$ be the measure preserving action of $\Gamma$ on $(Z,\rho)$ given by $$\gamma (x_1,x_2,\lambda)=(\gamma x_1,\gamma x_2,w(\gamma,x_1)\lambda w(\gamma,x_2)^{-1}),\forall (x_1,x_2)\in X\times X, \lambda\in\Lambda$$ and let $\pi:\Gamma\rightarrow \Cal U(L^2(Z,\rho))$ be the induced unitary representation. For every $n\geq 0$ and $a\in X_n$, let $\zeta_{a,n}\in L^{\infty}(Z,\rho)$ be the characteristic function of the set $X_{a,n}\times X_{a,n}\times\{e\}$ and define $\xi_{a,n}=\sqrt{|X_n|}\zeta_{a,n}$. Finally, for all $n\geq 0$, set $\xi_n=\sum_{a\in X_n}\xi_{a,n}$, then $||\xi_n||_{L^2(Z,\rho)}=1$. 
\vskip 0.1in
{\it Part 1.} $\lim_{n\rightarrow\infty}||\pi(\gamma)(\xi_n)-\xi_n||_{L^2(Z,\rho)}=0$, for all $\gamma\in\Gamma$.
\vskip 0.05in
{\it Proof of Part 1.} Fix $\gamma\in\Gamma$. For $\lambda\in\Lambda$, denote $S_{\lambda}=\{x\in X|w(\gamma^{-1},x)=\lambda\}$.  Using the fact that the set $\cup_{a\in X_n}(X_{a,n}\times X_{a,n})$ is invariant under the diagonal action $\alpha\times\alpha$ of $\Gamma$ on $X\times X$, we deduce that for all  $n$ $$<\pi(\gamma)(\xi_n),\xi_n>_{L^2(Z,\rho)}= \tag 3.a.$$

$$|X_n|\sum_{a\in X_n}\int_{X_{a,n}\times X_{a,n}\times \Lambda}\delta_{w(\gamma^{-1},x_1)\lambda w(\gamma^{-1},x_2)^{-1},e}\delta_{\lambda,e} d(\mu\times\mu\times c)(x_1,x_2,\lambda)=$$ $$|X_n|\sum_{a\in X_n}(\mu\times\mu)(\{(x_1,x_2)\in X_{a,n}\times X_{a,n}|w(\gamma^{-1},x_1)=w(\gamma^{-1},x_2)\})=$$ $$ |X_n|\sum_{a\in X_n,\lambda\in\Lambda}\mu(X_{a,n}\cap S_{\lambda})^2.$$
Now, let $P_n$ be the orthogonal projection from $L^2X$ onto the finite dimensional Hilbert space spanned by $\{1_{X_{a,n}}|a\in X_n\}.$ Then it is easy to check that for every measurable set $S\subset X$ and every $n\geq 0$, we have that $|X_n|\sum_{a\in X_n}\mu(X_{a,n}\cap S)^2=||P_n(1_S)||^2_{L^2(X,\mu)}$. By combining this fact with equality (3.a.)  we get that $$<\pi(\gamma)(\xi_n),\xi_n>_{L^2(Z,\rho)}=\sum_{\lambda\in\Lambda}||P_n(1_{S_{\lambda}})||_{L^2(X,\mu)}^2,\forall n\geq 0\tag 3.b.$$

Next, note that if we view $L^2X_n$ as a Hilbert subspace of $L^2X$ (via $r_n$), then $P_n$ is precisely the orthogonal projection onto $L^2X_n$. Since $X=\varprojlim X_n$ we have that $L^2X=\overline{(\cup_{n\geq 0}L^2X_n)}^{||.||}$, thus $P_n\rightarrow$ I, in the strong operator topology.

Now, let $\varepsilon>0$, then since $\sum_{\lambda\in\Lambda}\mu(S_{\lambda})=1$, we can find $F\subset\Lambda$ finite such that $\sum_{\lambda\in\Lambda\setminus F}\mu(S_{\lambda})\leq\varepsilon.$
Since $P_n\rightarrow$ I, (3.b.) implies that $$\liminf_{n} |<\pi(\gamma)(\xi_n),\xi_n>_{L^2(Z,\rho)}|\geq \sum_{\lambda\in F}||1_{S_{\lambda}}||_{L^2(X,\mu)}^2=\sum_{\lambda\in F}\mu(S_{\lambda})\geq 1-\varepsilon.$$
 As $\varepsilon>0$ is arbitrary and $||\xi_n||=1$, we get that $\lim_{n\rightarrow\infty}<\pi(\gamma)(\xi_n),\xi_n>_{L^2(Z,\rho)}=1$, for all $\gamma\in\Gamma$ and  Step 1 follows.\hfill$\square$
\vskip 0.03in
For the next three steps, fix $\delta\in (0,1)$. 

\vskip 0.1in
{\it Part 2.} There exists $n$  and $\eta\in L^2(Z,\rho)$ such that $\pi(\gamma)(\eta)=\eta$, for all $\gamma\in\Gamma_0$, and $||\eta-\xi_n||_{L^2(Z,\rho)}<\delta/4.$
\vskip 0.05in
{\it Proof of Part 2.}  Since the inclusion $\Gamma_0\subset\Gamma$ has relative property (T) (as defined in the introduction), then by [Jo1] we can find a finite set $F\subset\Gamma$ and $k>0$ such that if $\pi:\Gamma\rightarrow\Cal U(\Cal H)$ is a unitary representation and $\xi\in\Cal H$ is a unit vector with $||\pi(\gamma)(\xi)-\xi||\leq k$, for all $\gamma\in F$, then we can find a $\pi(\Gamma_0)$-invariant vector $\eta\in\Cal H$ such that $||\eta-\xi||<\delta/4$. Thus, using Part 1, we get the conclusion. 
\hfill$\square$

\vskip 0.1in
{\it Part 3. } There exists $a\in X_n$ such that $$||\pi(\gamma)(\xi_{a,n})-\xi_{a,n}||_{L^2(Z,\rho)}< (\delta/2)||\xi_{a,n}||_{L^2(Z,\rho)},\forall\gamma\in\Gamma_{a,n}\cap \Gamma_0.$$
\vskip 0.05in
{\it Proof of Part 3.} For every $a\in X_n$, let $\eta_a$ be the restriction of $\eta$ to $X_{a,n}\times X_{a,n}\times\Lambda$. 
Then $\eta'=\sum_{a\in X_n}\eta_a$ is the restriction of $\eta$ to $[\cup_{a\in X_n}(X_{a,n}\times X_{a,n})]\times\Lambda$.
Since $\xi_n$ is supported on $[\cup_{a\in X_n}(X_{a,n}\times X_{a,n})]\times\Lambda$, we get that $||\eta'-\xi_n||_{L^2(Z,\rho)}\leq ||\eta-\xi_n||_{L^2(Z,\rho)}< \delta/4$. Thus $$\sum_{a\in X_n}||\eta_a-\xi_{a,n}||_{L^2(Z,\rho)}^2=||\eta'-\xi_n||_{L^2(Z,\rho)}^2 < \delta^2/16=$$ $$(\delta^2/16)||\xi_{a,n}||_{L^2(Z,\rho)}^2=(\delta^2/16)\sum_{a\in X_n}  ||\xi_{a,n}||_{L^2(Z,\rho)}^2,$$ hence we can find $a\in X_n$ such that 
$||\eta_a-\xi_{a,n}||_{L^2(Z,\rho)} < (\delta/4) ||\xi_{a,n}||_{L^2(Z,\rho)}.$ 

On the other hand, since $\eta$ is $\pi(\Gamma_0)$-invariant and since  $X_{a,n}\times X_{a,n}\times\Lambda$ is a $\sigma(\Gamma_{a,n})$-invariant set, we get that $\eta_a$ is a $\pi(\Gamma_{a,n}\cap \Gamma_0)$-invariant vector. Thus by
 applying triangle's inequality we deduce that $$||\pi(\gamma)(\xi_{a,n})-\xi_{a,n}||_{L^2(Z,\rho)}\leq 2||\eta_{a}-\xi_{a,n}||_{L^2(Z,\rho)}<$$ $$(\delta/2)||\xi_{a,n}||_{L^2(Z,\rho)},\forall \gamma\in\Gamma_{a,n}\cap \Gamma_0.$$ \hfill$\square$
\vskip 0.1in
{\it Part 4.} There exists a  homomorphism $\psi:\Gamma_{a,n}\cap\Gamma_0\rightarrow\Lambda$ such that $w$ is cohomologous to a cocycle $w_1:\Gamma\times X\rightarrow \Lambda$ which satisfies $w_1(\gamma,x)=\psi(\gamma),$ for all $\gamma\in\Gamma_{a,n}\cap \Gamma_0$ and a.e. $x\in X_{a,n}$.
\vskip 0.05in
{\it  Proof of Part 4.} 
By  Part 3 we have that $$(1-\delta^2/8)||\xi_{a,n}||_{L^2(Z,\rho)}^2 <\Re <\pi(\gamma)(\xi_{a,n}),\xi_{a,n}>_{L^2(Z,\rho)},\forall\gamma\in\Gamma_{a,n}\cap \Gamma_0.\tag 3.c.$$

Now, a  computation similar to the one in Part 1  shows that inequality (3.c.) is equivalent to

$$(1-\delta^2/8)|X_n|^{-1}< |X_n|(\mu\times\mu)(\{(x_1,x_2)\in X_{a,n}\times X_{a,n}|w(\gamma,x_1)=w(\gamma,x_2)\}),$$ for all $\gamma\in\Gamma_{a,n}\cap \Gamma_0$. Hence, if we denote $\mu_{a,n}=(1/|X_n|)\mu_{|X_{a,n}}$, then $\mu_{a,n}$ is a probability measure on $X_{a,n}$ and the last inequality rewrites $$(\mu_{a,n}\times\mu_{a,n})(\{(x_1,x_2)\in X_{a,n}\times X_{a,n}|w(\gamma,x_1)=w(\gamma,x_2)\})>1-\delta^2/8,\forall\gamma\in\Gamma_{a,n}\cap \Gamma_0.$$
Next, since the action $\Gamma_0\curvearrowright (X,\mu)$ is assumed ergodic then by Remark 1.3.(3) we get that the action $\Gamma_{a,n}\cap\Gamma_0\curvearrowright (X_{a,n},\mu_{a,n})$ is ergodic. Thus, since we also have that $1-\delta^2/8>7/8$, we can apply Corollary 2.2. to deduce that there exists a group homomorphism $\psi:\Gamma_{a,n}\cap\Gamma_0\rightarrow\Lambda$ and a measurable function $\phi:X_{a,n}\rightarrow\Lambda$ such that  $\phi(\gamma x)w(\gamma,x)\phi(x)^{-1}=\psi(\gamma),$ for all $\gamma\in\Gamma_{a,n}\cap\Gamma_0$ and a.e. $x\in X_{a,n}.$ If we let $\tilde\phi:X\rightarrow\Lambda$ be given by $\tilde\phi(x)=\phi(x)$ if $x\in X_{a,n}$ and $\phi(x)=e$ otherwise, then $w_1(\gamma,x)=\tilde\phi(\gamma x)w(\gamma,x)\tilde\phi(x)^{-1}$ is a cocycle satisfying the conclusion.\hfill$\square$
\vskip 0.1in

{\it Part 5.} There exist homomorphisms $\psi_b:\Gamma_{b,n}\cap\Gamma_0\rightarrow\Lambda$, for all $b\in  X_n$, such that $w$ is cohomologous to a cocycle $w_2:\Gamma\times X\rightarrow \Lambda$ which satisfies $w_2(\gamma,x)=\psi_b(\gamma)$, for all $\gamma\in\Gamma_{b,n}\cap\Gamma_0$, $x\in X_{b,n}$ and $b\in X_n$. 
 
\vskip 0.05in
{\it Proof of Part 5.} Fix  $b\in X_n$ and  let $\gamma\in\Gamma$ such that  $b=\gamma a$. By using the cocycle relation we have that

$$w_1(\gamma_0,x)=w_1({\gamma}^{-1},\gamma_0 x)^{-1}w_1({\gamma}^{-1}\gamma_0\gamma,{\gamma}^{-1}x)w_1({\gamma}^{-1},x)\tag 3.d.$$ for all $\gamma_0\in\Gamma_{b,n}\cap \Gamma_0$ and a.e. $x\in X_{b,n}$.
 Since ${\gamma}^{-1}(\Gamma_{b,n}\cap\Gamma_0)\gamma=\Gamma_{a,n}\cap\Gamma_0$ ($\Gamma_0$ is normal in $\Gamma$) and ${\gamma}^{-1}X_{b,n}=X_{a,n}$, Part 4 implies that $w_1({\gamma}^{-1}\gamma_0\gamma,{\gamma}^{-1}x)=\psi({\gamma}^{-1}\gamma_0\gamma),$ for all  $\gamma_0\in\Gamma_{b,n}\cap\Gamma_0$ and a.e. $x\in X_{b,n}.$ Define $\phi_b:X_{b,n}\rightarrow\Lambda$ by $\phi_b(x)=w_1({\gamma}^{-1},x)
$ and $\psi_b:\Gamma_{b,n}\cap\Gamma_0\rightarrow\Lambda$ by $\psi_b(\gamma_0)=\psi({\gamma}^{-1}\gamma_0\gamma)$. Then $\psi_b$ is a homomorphism and identity (3.d.) becomes $$w_1(\gamma_0,x)=\phi_b(\gamma_0 x)^{-1}\psi_b(\gamma)\phi_b(x),\forall\gamma_0\in\Gamma_{b,n}\cap\Gamma_0\tag 3.e.$$ and a.e. $x\in X_{b,n}$.

 Finally, let $\phi:X\rightarrow\Lambda$ be defined by $\phi(x)=\phi_b(x)$ iff $x\in X_{b,n}$. Then the formula $w_2(\gamma,x)=\phi(\gamma x)w_1(\gamma,x)\phi(x)^{-1}$ defines a cocycle cohomologous to $w_1$ (thus to $w$) which by (3.e.) verifies the conclusion.\hfill$\square$
\vskip 0.1in
{\it Part 6.} There exists $N\geq n$ such that $w_2$ factors through the map 
$\Gamma\times X\rightarrow \Gamma\times X_N$.
\vskip 0.05in
{\it Proof of Part 6.}
Fix $\gamma\in\Gamma$. We claim first that there exists $N(\gamma)\geq n$ such that map $X\ni x\rightarrow w_2(\gamma,x)\in\Lambda$ factors through $r_{N(\gamma)}:X\rightarrow X_{N(\gamma)}$.
Start by noticing that Part 5 and identity (3.d.) for $w_2$ imply that $$w_2(\gamma,\gamma_0 x)=\psi_{\gamma b}(\gamma\gamma_0{\gamma}^{-1})w_2(\gamma,x)\psi_b(\gamma_0)^{-1},$$for all  $b\in X_n,\gamma_0\in\Gamma_{b,n}\cap\Gamma_0$ and a.e. $x\in X_{b,n}.$ In particular, for all $b\in X_n$, the map $X_{b,n}\ni x\rightarrow w_2(\gamma,x)\in\Lambda$
satisfies the hypothesis of Corollary 1.5. for the action $(\Gamma_{b,n}\cap\Gamma_0)\curvearrowright (X_{b,n},\mu_{b,n}).$ 
Now, by Remark 1.3.(3) this action is ergodic and profinite (since the action $\Gamma_0\curvearrowright (X,\mu)$ is ergodic and profinite). Corollary 1.5. thus implies that for all $b\in X_n$, we can find $N(\gamma,b)\geq n$ such that the map $x\mapsto w_2(\gamma,x)$ is constant on $X_{b',N(\gamma,b)}$, for all $b'\in X_{N(\gamma,b)}$ such that $X_{b',N(\gamma,b)}\subset X_{b,n}$.  It is now clear that the claim holds true for $N(\gamma)=\max_{b\in X_n}N(\gamma,b).$

Next, if we let $\Gamma^n=\cap_{a\in X_n}\Gamma_{a,n}$, then Part 5 gives that the maps $x\rightarrow w_2(\gamma_0,x)$ are constant on $X_{b,n}$, for all $b\in X_n$ and all $\gamma_0\in\Gamma^n\cap\Gamma_0$.   Since $[\Gamma_0:(\Gamma^n\cap\Gamma_0)]<\infty$ and $\Gamma/\Gamma_0$ is finitely generated we can find   $\gamma_1,...,\gamma_m\in\Gamma$ such that $\Gamma$ is generated by $\Gamma^n\cap\Gamma_0$ and $\gamma_1,..,\gamma_m$. 
Let $N=\max_{i=\overline{1,m}}N(\gamma_i)$, then the maps $x\rightarrow w_2(\gamma,x)$ are constant on $X_{b,N}$, for all $b\in X_N$ and for all  $\gamma\in S:=(\Gamma_{a,n}\cap\Gamma_0)\cup\{\gamma_1,..,\gamma_m\}$. 

Finally, let $g\in \Gamma$, then since  $S$ generates $\Gamma$ we can find $g_1,..,g_l\in S$ such that $g=g_1g_2...g_l$. By the cocycle identity we have that a.e. $x\in X$ $$w(g,x)=w(g_1,(g_2...g_l)x)w(g_2,(g_3...g_l)x)...w(g_l,x)\tag 3.f.$$ Let $b\in X_N$ and $2\leq k\leq l$. Then $(g_kg_{k+1}..g_l)X_{b,n}=X_{(g_kg_{k+1}..g_l)b,n}$ and since $g_{k-1}\in S$, the above implies that  the map $X_{b,N}\ni x\rightarrow w_2(g_{k-1},(g_kg_{k+1}..g_l)x)$ is constant. Combining this fact with identity (3.f.) we get that the map $X_{b,N}\ni x\rightarrow w_2(g,x)$ is constant on $X_{b,N}$, for all $b\in X_N$.  \hfill$\blacksquare$
\vskip 0.05in
\noindent
{\bf 3.1. Remarks.} (1). Note that in the above proof we do not use the full strength of the relative property (T) of the inclusion $\Gamma_0\subset\Gamma$. In fact, notice that  in order to untwist $w$ we only use relative property (T) for the specific representation $\pi$ arising from the action $\Gamma\curvearrowright^{\sigma} X\times X\times\Lambda$. 

(2). We remark that  Theorem B in combination with Popa's cocycle superrigidity result ([P2]) leads to examples of cocycle superrigid actions which are neither profinite (as in Theorem B) nor Bernoulli (as in [P2]). For this, let $\Gamma$, $\Gamma_0$ and $\Gamma\curvearrowright^{\alpha}X$ be as in Theorem A. Let $\Gamma\curvearrowright^{\rho}Y=[0,1]^{\Gamma}$ be the Bernoulli shift action (where $[0,1]$ is endowed with the Lebesgue measure). Finally, let $w:\Gamma\times (X\times Y)\rightarrow \Lambda$ be a cocycle for the diagonal product action $\alpha\times\rho$ with values in a countable group $\Lambda$. Then there exist $n$ such that $w$ is cohomologous to a cocycle $w':\Gamma\times (X\times Y)\rightarrow\Lambda$ of the form $w'=w''\circ(id\times s_n)$, for some cocycle $w'':\Gamma\times X_n\rightarrow\Lambda$, where $s_n:X\times Y\rightarrow X_n$ is given by $s_n(x,y)=r_n(x)$. Indeed, by Theorem 5.2. in [P2], $w$ is cohomologous to a cocycle which only depends on the $X$-variable and then Theorem B gives the claimed conclusion.

(3). Recently, S. Coskey  used Theorem B to prove that the isomorphism and quasi-isomorphism
problems for the p-local torsion-free abelian groups of rank n have incomparable Borel complexities, for every $n\geq 3$ ([Co]).

\vskip 0.1in
\head \S 4. Proof of OE superrigidity.\endhead

In this Section we prove a more general version of Theorem A. To state it, we need to review some more terminology (see, for example, Section 1 in [P2]).
Let $\Cal R$ be a countable measure preserving equivalence relation on a standard probability space $(X,\mu)$. Recall that every such equivalence relation is of the form $\Cal R_{\Gamma}=\{(x,\gamma x)|x\in X,\gamma\in\Gamma\}$, for some measure preserving action $\Gamma\curvearrowright X$ of a countable group $\Gamma$ ([FM]). Hereafter, we will refer to the equivalence relation $\Cal R_{\Gamma}$ as the equivalence relation induced by the action $\Gamma\curvearrowright X$. 

Now,  assume that $\Cal R$ is an ergodic equivalence relation and let $t>0$. On $\tilde X=X\times\Bbb N$, endowed with the measure
$\tilde\mu=\mu\times c$ (where $c$ is the counting measure on $\Bbb N$), consider the equivalence relation $\tilde{\Cal R}$ given by $(x,i)\tilde{\Cal R}(y,j)$ $\Longleftrightarrow$ $x\Cal R y$. Let $Y\subset \tilde X$ be a set of measure $\tilde\mu(Y)=t$ and define $\Cal R^Y=\tilde{\Cal R}\cap (Y\times Y)$ (the restriction of $\tilde{\Cal R}$ to $Y$). Then $\Cal R^Y$ is a countable equivalence relation on $Y$ preserving the probability measure $t^{-1}\tilde{\mu}_{|Y}$. Moreover, its isomorphism class, denoted $\Cal R^t$ (in words: the {\bf t-amplification of $\Cal R$}), only depends on $t$ and not on the choice of $Y$.

Also,  we say that  $\theta:\Cal R\rightarrow \Cal S$ is an {\bf orbit equivalence} (or {\bf isomorphism}) between   two countable measure preserving equivalence relations  $\Cal R,\Cal S$ on  $(X,\mu)$ and on $(Y,\nu)$, respectively, if $\theta:X\rightarrow Y$ is a probability space isomorphism such that $\theta$ is a bijection between the $\Cal R$-orbit of $x$ and the $\Cal S$-orbit of $\theta(x)$, a.e. $x\in X$.  
Finally, we recall that the {\bf full group} of  a countable measure preserving equivalence relation $\Cal R$ on $(X,\mu)$, denoted $[\Cal R]$, consists of the automorphisms $\tau$ of $X$ such that $\tau(x)\Cal R x$, a.e. $x\in X$.
 
\proclaim {4.1. Theorem}    Let $\Gamma\curvearrowright^{\alpha}X$ be as in Theorem A.  Suppose that $\alpha$ is the limit of the actions $\Gamma\curvearrowright^{\alpha_n} X_n$, with $X_n$ finite, and, for every $n$, let $r_n:X\rightarrow X_n$ be the quotient map. For every $n$ and $a\in X_n$, denote $X_{a,n}=r_n^{-1}(\{a\})$ and $\Gamma_{a,n}=\{\gamma\in\Gamma|\gamma X_{a,n}=X_{a,n}\}$.

Let $\Lambda\curvearrowright^{\beta} Y$ be a free ergodic measure preserving action of a countable group $\Lambda$ on a standard probability space $Y$ and let $\theta:\Cal R_{\Gamma}\rightarrow \Cal R_{\Lambda}^t$ be an orbit equivalence, for some $t>0$.    
Then we can find $n$, $a\in X_n$, a finite index subgroup $\Lambda_0$ of $\Lambda$, $\tau\in[\Cal R_{\Lambda}]$ and a group isomorphism $\psi:\Gamma_{a,n}\rightarrow \Lambda_0$ such that $Y_0=(\tau\circ\theta)(X_{a,n})$ is a $\Lambda_0$-invariant set and $(\tau\circ\theta)_{|X_{a,n}}:X_{a,n}\rightarrow  Y_0$ conjugates the actions  $\Gamma_{a,n}\curvearrowright X_{a,n}$ and $\Lambda_0\curvearrowright Y_0$. 
Moreover, $\mu(Y_0)=[\Lambda:\Lambda_0]^{-1}$, $t=[\Lambda:\Lambda_0]^{-1}[\Gamma:\Gamma_{a,n}]\in\Bbb Q$ and the action $\Lambda\curvearrowright Y$ is obtained by inducing the action $\Lambda_0\curvearrowright Y_0$ to $\Lambda$.
\endproclaim

{\it Proof.}
We first assume that $t\leq 1$. Let $Y'\subset Y$ be a measurable set of measure $t$ and let $\theta:(X,\mu)\rightarrow (Y',t^{-1}\nu_{|Y'})$ be a probability space isomorphism such that $\theta(\Gamma x)=\Lambda\theta(x)\cap Y',$ a.e. $x\in X$.
For $x\in X$ and $\gamma\in\Gamma$, let $w(\gamma,x)$ be the unique (by freeness of the $\Lambda$-action) element of $\Lambda$ such that  $\theta(\gamma x)=w(\gamma,x)\theta(x)$. Then $w:\Gamma\times X\rightarrow\Lambda$ is a measurable cocycle. By applying Theorem B we can find $n\geq 0$, $a\in X_n$, a group homomorphism $\psi:\Gamma_{a,n}\rightarrow\Lambda$ and a measurable map $\phi:X\rightarrow\Lambda$ such that   $$w(\gamma,x)=\phi(\gamma x)\psi(\gamma)\phi(x)^{-1}\tag 4.a.$$ for all  $\gamma\in\Gamma_{a,n}$ and a.e. $x\in X_{a,n}$. 

\vskip 0.05in
Next, let $\lambda\in\Lambda$ such that $\mu(\{x\in X_{a,n}|\phi(x)=\lambda\})>0$. Since $X=\varprojlim X_m$, we can find $m\geq n$ and $b\in X_m$ such that $X_{b,m}\subset X_{a,n}$ and $\mu(\{x\in X_{b,m}|\phi(x)=\lambda\})>(1/2)\mu(X_{b,m})$.
 Define $\phi'(x)=\phi(x)\lambda^{-1}$ and $\psi'(\gamma)=\lambda\psi(\gamma){\lambda}^{-1}$ for all $x\in X$ and $\gamma\in\Gamma_{a,n}$.
 Thus, if we replace $n,a,\phi,\psi$ by $m,b,\phi',\psi'$, respectively, then identity (4.a.) still holds true and we moreover have that $\mu(\{x\in X_{a,n}|\phi(x)=e\})>(1/2)\mu(X_{a,n})$.
  
\vskip 0.05in
 Let $n$ and $a\in X_n$ be as above. To simplify notation, set $A=\{x\in X_{a,n}|\phi(x)=e\}$, $B=\theta(A)$,  and $\Lambda_0=\psi(\Gamma_{a,n})$. 
 Also,  define $\tilde\theta:X\rightarrow Y$ by $\tilde\theta(x)=\phi(x)^{-1}\theta(x)$, for all $x\in X$, and let $Y_0=\tilde\theta(X_{a,n})$.
Then  from (4.a.) and the above discussion we derive the following two relations $$\tilde\theta(\gamma x)=\phi(\gamma x)^{-1}\theta(\gamma x)=\phi(\gamma x)^{-1}w(\gamma,x)\theta(x)=\tag 4.b.$$ $$\psi(\gamma)\phi(x)^{-1}\theta(x)=\psi(\gamma)\tilde\theta(x), \forall \gamma\in\Gamma_{a,n}, $$ a.e. $x\in X_{a,n}$ and $$\mu(A)>(1/2)\mu(X_{a,n})\tag 4.c.$$  

\vskip 0.1in

{\it Claim 1.} $\tilde\theta_{|X_{a,n}}$ is 1-1, $\nu(Y_0)=t\mu(X_{a,n})$,  $$\tilde\theta_{|X_{a,n}}:(X_{a,n},{\mu(X_{a,n})}^{-1}\mu_{|X_{a,n}})\rightarrow(Y_0,{\nu(Y_0)}^{-1}\nu_{|Y_0})$$ is a probability space isomorphism and Ker$(\psi)=\{e\}$.

\vskip 0.05in

{\it Proof of Claim 1.} We first note that if $C\subset X$ is a measurable set such that $\tilde\theta_{|C}$ is 1-1, then $\nu(\tilde\theta(C))=t\mu(C)$. Indeed, if $C_{\lambda}=\{x\in C|\phi(x)=\lambda\}$, then $C=\sqcup_{\lambda\in\Lambda}C_{\lambda}$ and $\tilde\theta(C)=\sqcup_{\lambda\in\Lambda}\lambda\theta(C_{\lambda})$, thus $$\nu(\tilde\theta(C))=\sum_{\lambda\in\Lambda}\nu(\lambda\theta(C_{\lambda}))=\sum_{\lambda\in\Lambda}\nu(\theta(C_\lambda))=t\sum_{\lambda\in\Lambda}\mu(C_\lambda)=t\mu(C).$$

Next, since $B=\theta(A)=\tilde\theta(A)\subset Y_0,$ we get that $\nu(Y_0)\geq \nu(B)$. Since $\nu(B)=t\mu(A)$, relation (4.c.)  implies that $\nu(Y_0)>(t/2)\mu(X_{a,n}).$
On the other hand,  (4.b.) implies that the function $$X_{a,n}\ni x\rightarrow |\{y\in X_{a,n}|\tilde\theta(x)=\tilde\theta(y)\}|\in\Bbb N\cup\{+\infty\}$$ is $\Gamma_{a,n}$-invariant. Since the action $\Gamma_{a,n}\curvearrowright (X_{a,n},\mu_{a,n})$ is ergodic (by Remark 1.3.(3)), we can find $k\geq 1$ such that $|\{y\in X_{a,n}|\tilde\theta(x)=\tilde\theta(y)\}|=k$, a.e. $x\in X_{a,n}$.

If $k\geq 2$, then we can find two disjoint measurable sets $Z_1,Z_2\subset X_{a,n}$  such that $\tilde\theta(Z_i)=Y_0$ and 
$\tilde\theta_{|Z_i}$ is 1-1, for $i\in\{1,2\}$.
Indeed, note that on one hand the above gives that $|{\tilde{\theta}}_{|X_{a,n}}^{-1}(\{y\})|\geq 2,$ 
a.e. $y\in Y_0$, while, on the other hand ${\tilde\theta}^{-1}(\{y\})\subset 
\Gamma\theta^{-1}(y)$, a.e. $y\in Y$ (here we use the definition of
$\tilde\theta$ and the fact that $\theta(\Gamma x)=\Lambda\theta(x)\cap Y'$, a.e. $x\in X$).
Thus, a.e. $y\in Y_0$, we can find $\gamma_{1,y},\gamma_{2,y}\in\Gamma$ (depending measurably on $y$)
 such that
$\gamma_{1,y}\theta^{-1}(y),\gamma_{2,y}\theta^{-1}(y)\in 
{\tilde{\theta}}_{|X_{a,n}}^{-1}(\{y\})$. In this then clear 
that $Z_i=\{\gamma_{i,y}\theta^{-1}(y)|y\in Y_0\}$, $i=1,2,$ verify the desired conditions.
 However, by the first part of the proof of this claim we would get 
 that $\mu(Z_1)=\mu(Z_2)=(1/t)\nu(Y_0)$, thus $$\mu(X_{a,n})\geq\mu(Z_1)+\mu(Z_2)=(2/t)\nu(Y_0),$$ a contradiction. Therefore $k$ must be 1, hence  $\tilde\theta_{|X_{a,n}}$ is 1-1.

  For the last part, let $\gamma\in$ Ker$(\psi)$. Then by identity (4.b.) we get that $\tilde\theta(\gamma x)=\tilde\theta(x)$ a.e. $x\in X_{a,n}$. Since $\tilde\theta$ is 1-1 on $X_{a,n}$, it follows that $\gamma x=x$, a.e. $x\in X_{a,n}$ and since $\Gamma$ acts freely on $X$, we must have $\gamma=e$.  \hfill$\square$

\vskip 0.1in
{\it Claim 2.} $\Lambda_0$ is a finite index subgroup of $\Lambda$.
\vskip 0.05in
{\it Proof of Claim 2.} To prove this, we use an argument from Section 5 in [Fu4]. Let $\lambda\in\Lambda$ such that $\nu(B\cap \lambda^{-1} B)>0$. We claim that $\lambda\in\Lambda_0$. 
To this end, let $y\in B\cap \lambda^{-1}B$. Then we can find $x_1,x_2\in A\subset X_{a,n}$ such that $\theta(x_1)=y$ and $\theta(x_2)=\lambda y$, thus $\theta(x_2)\in\Lambda\theta(x_1)\cap B\subset \Lambda\theta(x_1)\cap\theta(X_{a,n})$. Using the hypothesis we have that $$\theta(\Gamma_{a,n}x)=\theta(\Gamma x\cap X_{a,n})=\theta(\Gamma x)\cap\theta(X_{a,n})=\Lambda\theta(x)\cap\theta(X_{a,n}),$$ a.e. $x\in X_{a,n}$. Thus, we can find $\gamma\in\Gamma_{a,n}$ such that $\theta(x_2)=\theta(\gamma x_1)$ and since $\theta$ is 1-1, we derive that $x_2=\gamma x_1$. Further, using (4.b.), we get that $$\theta(x_2)=\tilde\theta(x_2)=\tilde\theta(\gamma x_1)=\psi(\gamma)\tilde\theta(x_1)=\psi(\gamma)\theta(x_1)$$ and since we also have that $\theta(x_2)=\lambda\theta(x_1)$ freeness of the action of $\Lambda$ on $Y$ implies that $\lambda=\psi(\gamma)\in\Lambda_0$. 

Finally, if we assume that $[\Lambda:\Lambda_0]=\infty$, then we can find $\lambda_1,\lambda_2,...\in\Lambda$ such that ${\lambda_i}^{-1}\lambda_j\notin\Lambda_0$ for all $i\not=j$. 
The above then implies that $\nu(\lambda_i B\cap \lambda_j B)=\nu(B\cap {\lambda_i}^{-1}\lambda_j B)=0$ whenever $i\not=j$. Thus, since $\lambda_i B\subset Y$, for all $i$, we get that $\mu(Y)\geq\sum_{i\geq 1}\nu(\lambda_i B)=\sum_{i\geq 1}\nu(B)=+\infty$, a contradiction.\hfill $\square$
\vskip 0.1in 	

{\it Claim 3.}  $\Lambda_0y=\Lambda y\cap Y_0$, a.e. $y\in Y_0$.
\vskip 0.05in
{\it Proof of Claim 3.} By the hypothesis and the fact that $\theta(A)=B$, we deduce that $$\theta(\Gamma x\cap A)=\Lambda\theta(x)\cap B\tag 4.d.$$a.e. $x\in X_{a,n}$. On the other hand, if $x,\gamma x\in A\subset X_{a,n}$, for some $\gamma\in\Gamma$, then $\gamma\in\Gamma_{a,n}$ and $\theta(\gamma x)=\tilde\theta(\gamma x)=\psi(\gamma)\tilde\theta(x)=\psi(\gamma)\theta(x),$ thus $\theta(\gamma x)\in\Lambda_0\theta(x)$. Altogether, we get that $$\theta(\Gamma x\cap A)\subset\Lambda_0\theta(x)\cap B\tag 4.e.$$ a.e. $x\in A$. Combining (4.d.) and (4.e.) we get that $ \Lambda_0\theta(x)\cap B\supset\Lambda\theta(x)\cap B$, a.e. $x\in A$, thus $$ \Lambda_0y\cap B=\Lambda y\cap B\tag 4.f.$$a.e. $y\in B$. Next, note that by (4.b.) and Claim 1, the actions $$\Gamma_{a,n}\curvearrowright (X_{a,n},\mu_{|X_{a,n}}/\mu(X_{a,n})),\Lambda_0\curvearrowright (Y_0,\nu_{|Y_0}/\nu(Y_0)$$ are conjugate. In particular, the action of $\Lambda_0$ on $Y_0$ is ergodic.

Now, let $y\in Y_0$. Since $\Lambda_0 y\subset\Lambda y\cap Y_0$ ($Y_0$ is $\Lambda_0$-invariant), to prove the claim we only need to show that if $\lambda\in\Lambda$ is such that $\lambda y\in Y_0$, then $\lambda\in\Lambda_0$.
 Since $B\subset Y_0$ is of positive measure and $\Lambda_0$ acts ergodically on $Y_0$ we can find $\lambda_1,\lambda_2\in\Lambda_0$
such that $\lambda_1y,\lambda_2(\lambda y)\in B$. This implies that $(\lambda_2\lambda)y\in\Lambda(\lambda_1y)\cap B$ and since $\lambda_1 y\in B$, relation (4.f.) gives that we can find $\lambda_3\in\Lambda_0$ such that $(\lambda_2\lambda)y=\lambda_3(\lambda_1 y)$. Finally, the freeness of the action of $\Lambda$ on $Y$ implies that $\lambda_2\lambda=\lambda_3\lambda_1$, hence $\lambda={\lambda_2}^{-1}\lambda_3\lambda_1\in\Lambda_0.$\hfill $\square$	
\vskip 0.1in

We can now finalize the proof in the case $t\leq 1$. 
Firstly, observe that by Claim 1, $\tau'=\tilde\theta\circ{\theta}^{-1}: (\theta(X_{a,n}),\nu_{|\theta(X_{a,n})})\rightarrow (Y_0,\nu_{|Y_0})$ is a probability space isomorphism with the property that $\tau'(y)\in\Lambda y$, a.e. $y\in \theta(X_{a,n})$. Since $\Lambda$ acts ergodically on $Y$, we can find $\tau\in [\Cal R_{\beta}]$ such that $\tau'=\tau_{|\theta(X_{a,n})}$. Then by formula (4.b.) we get that $Y_0=(\tau\circ\theta)(X_{a,n})$ is a $\Lambda_0$-invariant set and that ${\tau\circ\theta}_{|X_{a,n}}$ conjugates the actions $\Gamma_{a,n}\curvearrowright X_{a,n}$ and $\Lambda_0\curvearrowright Y_0$.

Secondly, notice that since $\Lambda$ acts freely on $Y$, then by Claim 3 we get that $\nu(\lambda Y_0\cap Y_0)=0$ if $\lambda\not\in\Lambda_0$. On the other hand, since the action $\Lambda\curvearrowright Y$ is ergodic, we get that $Y=\cup_{\lambda\in\Lambda}\lambda Y_0$, a.e. By combining these two observations, it follows that if $\lambda_1,..,\lambda_k\in\Lambda$ are such that $\Lambda=\sqcup_{i=1}^k\lambda_i\Lambda_0$, then $Y$ is the disjoint union of $\lambda_1 Y_0,..,\lambda_k Y_0$. Thus $\nu(Y)=\sum_{i=1}^k\nu(\lambda_i Y_0)=[\Lambda:\Lambda_0]\nu(Y_0)$, hence $\nu(Y_0)=[\Lambda:\Lambda_0]^{-1}$. Also, by Claim 1, $\nu(Y_0)=t\mu(X_{a,n})=t[\Gamma:\Gamma_{a,n}]^{-1}$. Altogether, we get that $t=[\Lambda:\Lambda_0]^{-1}[\Gamma:\Gamma_{a,n}]$.

Using again the fact that   $Y$ is the disjoint union of $\lambda_1 Y_0,..,\lambda_k Y_0$, it is easy to see that the action $\Lambda\curvearrowright Y$ is obtained by inducing the action $\Lambda_0\curvearrowright Y_0$ to $\Lambda$.

\vskip 0.05in
In the case $t>1$, we first remark that if $\theta:X\rightarrow Y$ is an orbit equivalence between two  equivalence relations $\Cal R$ and $\Cal S$, then the restriction of $\theta$ to a set $X_0\subset X$ of measure $s$ induces an orbit equivalence between $\Cal R\cap (X_0\times X_0)$ and $\Cal S^s$. Thus, if $n\geq 0$ is such that $|X_n|\geq t$, then the restriction of $\theta$ to $X_{a,n}$ (for some $a\in X_n$) induces an orbit equivalence between $\Cal R_{\Gamma}\cap (X_{a,n}\times X_{a,n})$ and $\Cal R_{\Lambda}^{t|X_n|^{-1}}$. On the other hand,  $\Cal R_{\Gamma}\cap (X_{a,n}
\times X_{a,n})$ is precisely the equivalence relation induced by the free ergodic profinite action $\Gamma_{a,n}\curvearrowright X_{a,n}$. Note that this action trivially verifies the hypothesis of Theorem 4.1.
Since we also have that $t|X_n|^{-1}\leq 1$, the first part of the proof gives the conclusion in this case. 
\hfill $\blacksquare$

\head \S 5. Applications.\endhead
\vskip 0.1in
We derive in this Section several consequences of Theorem 4.1. To start with, we give a construction of a concrete uncountable family of non-orbit equivalent profinite actions for SL$_n(\Bbb Z)$ ($n\geq 3$) as well as for its finite index subgroups.
To this end, fix $n\geq 3$ and  a finite index subgroup $\Gamma\subset$ SL$_n(\Bbb Z)$ and remark that since SL$_n(\Bbb Z)$ has property (T), $\Gamma$ also has property (T) ([K]). For every $m\geq 1$,  denote SL$_n(m\Bbb Z)=\ker($SL$_n(\Bbb Z)\rightarrow$ SL$_n(\Bbb Z/m\Bbb Z)).$ Let $I=\{p_k\}_{k\geq 1}$ be an infinite increasing sequence of prime numbers. For all $k\geq 1$,  define $$\Gamma_{I,k}=\Gamma\cap {\text {SL}}_n({p_1p_2..p_k}\Bbb Z)$$   and let $\alpha_I$ be the measure preserving profinite action $\Gamma\curvearrowright (X_I,\mu_I):=\varprojlim(\Gamma/\Gamma_{I,k},\mu_k)$, where $\mu_k$ is the  counting probability measure on $\Gamma/\Gamma_{I,k}$. Since $\cap_{k\geq 1}$ SL$_n(p_1p_2..p_k\Bbb Z)=\{1\}$, we have that
 $\alpha_I$ is  free and ergodic.
Next, we show that (essentially) different sets of primes $I_1$ and $I_2$ give to non conjugate (not even "virtually" conjugate)  actions $\alpha_{I_1}$ and $\alpha_{I_2}$. When combined with Theorem 4.1., this implies that the actions $\alpha_{I_1}$ and $\alpha_{I_2}$ are not stably orbit equivalent.

\proclaim {5.1. Proposition} Let $I_1=\{p_k^1\}_{k\geq 1}$ and $I_2=\{p_k^2\}_{k\geq 1}$ be two infinite sequences of primes and let $\Gamma_1,\Gamma_2\subset \Gamma$ be two finite index subgroups. For $j\in\{1,2\}$, let $(Y_j,\nu_j)\subset (X_{I_j},\mu_{I_j})$ be an ergodic component for the restriction of  ${\alpha_{I_j}}$ to ${\Gamma_j}$. If the actions $\Gamma_1\curvearrowright (Y_1,\nu_1)$ and $\Gamma_2\curvearrowright (Y_2,\nu_2)$ are conjugate, then $|I_1\Delta I_2|<\infty$. 
\endproclaim
{\it Proof.} Assume that the actions $\Gamma_1\curvearrowright (Y_1,\nu_1)$ and $\Gamma_2\curvearrowright (Y_2,\nu_2)$ are conjugate. We will prove that $|I_1\setminus I_2|<\infty$. Similarly, it follows that $|I_2\setminus I_1|<\infty$, thus giving the conclusion. Note first that it is easy to see that for $j\in\{1,2\}$, the action $\Gamma_j\curvearrowright (Y_j,\nu_j)$ is isomorphic to the action $$\Gamma_j\curvearrowright\varprojlim \Gamma_j/[\Gamma_j\cap  {\text{SL}}_n(p_1^j..p_k^j\Bbb Z))].$$ 
Thus, by applying Lemma 1.8. we can find a group isomorphism $\delta:\Gamma_1\rightarrow\Gamma_2$ and two subsequences $\{n_k\}_k,\{N_k\}_k$ such that  $$\delta(\Gamma_1\cap {\text{SL}}_n(p_1^1..p_{n_{k}}^1\Bbb Z))\subset \Gamma_2\cap  {\text{SL}}_n(p_1^2..p_{N_k}^2\Bbb Z),\forall k\tag 5.1.a.$$
Next, we need the following lemma, whose proof, although standard,  we have included below for completeness. 
\proclaim {5.2. Lemma}
Let $n\geq 3$ and $\Gamma,\Lambda\subset{\text{SL}}_n(\Bbb Z)$ be two finite index subgroups. Assume that $\delta:\Gamma\rightarrow\Lambda$ is a group isomorphism. Then there exists an invertible matrix $A\in \Bbb M_n(\Bbb Z)$ such that either 
$(i)$ $\delta(\gamma)=A\gamma A^{-1},$ for all $\gamma\in\Gamma$, or 
$(ii)$ $\delta(\gamma)=A(\gamma^{-1})^{t}A^{-1},$ for all $\gamma\in\Gamma$.
\endproclaim

Going back to the proof of Proposition 5.1., apply Lemma 5.2. to the group isomorphism $\delta:\Gamma_1\rightarrow\Gamma_2$. Thus,  we can find $A\in \Bbb M_n(\Bbb Z)$ invertible such that, after eventually replacing $\Gamma_1$ with $\Gamma_1^t=\{\gamma^t|\gamma\in\Gamma_1\}$, we have that  for all $k$ $$A[\Gamma_1\cap {\text{SL}}_n(p_1^1..p_{n_{k}}^1\Bbb Z)]A^{-1}\subset \Gamma_2\cap {\text{SL}}_n(p_1^2..p_{N_k}^2\Bbb Z)\tag 5.1.b.$$ 
Let $l$ be the largest number such that  $p_{l}^2|\det A$. Then for all $k>l$ we have that  $$A^{-1}{\text{SL}}_n(p_1^2..p_{N_k}^2\Bbb Z)A\cap {\text{SL}}_n(\Bbb Z)\subset {\text{SL}}_n(p_{{l+1}}^2..p_{N_k}^2\Bbb Z)\tag 5.1.c.$$ 
By combining (5.1.b.) and (5.1.c.) we deduce that $$\Gamma_1\cap {\text{SL}}_n(p_1^1..p_{n_{k}}^1\Bbb Z)\subset {\text{SL}}_n(p_{{l+1}}^2..p_{N_k}^2\Bbb Z),\forall k>l.$$
This further gives that $$[{\text{SL}}_n(p_1^1..p_{n_{k}}^1\Bbb Z):({\text{SL}}_n(p_1^1..p_{n_{k}}^1\Bbb Z)\cap {\text{SL}}_n(p_{{l+1}}^2..p_{N_k}^2\Bbb Z))]\leq\tag 5.1.d.$$
$$[{\text{SL}}_n(p_1^1..p_{n_{k}}^1\Bbb Z):({\text{SL}}_n(p_1^1..p_{n_{k}}^1\Bbb Z)\cap \Gamma_1)]\leq [{\text{SL}}_n(\Bbb Z):\Gamma_1],\forall k>l.$$
Now, for every $k>l$, define $P_k=\{p_1^1,..,p_{n_{k}}^1\}\setminus\{p_{l+1}^2,..,p_{N_k}^2\}$.
Then it is easy to see that (5.1.d.) is equivalent to $$\prod_{p\in P_k}|{\text{SL}}_n(\Bbb Z/p\Bbb Z)|\leq [{\text{SL}}_n(\Bbb Z):\Gamma_1],\forall k>l\tag 5.1.e.$$Since $\lim_{p\rightarrow\infty}|{\text{SL}}_n(\Bbb Z/p\Bbb Z)|=+\infty$, by using (5.1.e.), we get that we can find $N\in\Bbb N$ such that $$\{p_1^1,..,p_{n_{k}}^1\}\setminus\{p_{l+1}^2,..,p_{N_k}^2\}\subset\{2,3,..,N\},\forall k>l.$$ Finally, since $I_1\setminus I_2\subset\cup_{k>l}(\{p_1^1,..,p_{n_{k}}^1\}\setminus\{p_{l+1}^2,..,p_{N_k}^2\})$, we get that $|I_1\setminus I_2|<\infty$.
 \hfill$\blacksquare$
\vskip 0.1in
{\it Proof of Lemma 5.2.} Since $\Gamma$ and $\Lambda$ are finite index subgroups of SL$_n(\Bbb Z)$, they are lattices in SL$_n(\Bbb R)$.  Also, since $n\geq 3$,  Mostow's rigidity theorem (see [Z], for example) implies that $\delta$ extends to an automorphism of SL$_n(\Bbb R)$. Thus, we can find $A\in$ SL$_n(\Bbb R)$ such that either $(i)$ $\delta(\gamma)=A\gamma A^{-1},$ for all $\gamma\in$ SL$_n(\Bbb R),$ or 
$(ii)$ $\delta(\gamma)=A(\gamma^{-1})^{t}A^{-1},$ for all $\gamma\in$ SL$_n(\Bbb R)$. Next, if we denote $\Gamma_0=\{\gamma\in\Gamma|\gamma^t\in\Gamma\}$, then $\Gamma_0\subset$ SL$_n(\Bbb Z)$ is a finite index subgroup and in any of the above cases we have that $A\Gamma_0 A^{-1}\subset$ SL$_n(\Bbb Z)$.

Since any finite index subgroup of SL$_n(\Bbb Z)$ ($n\geq 3$) is a congruence subgroup, we can find $m\geq 1$ such that SL$_n(m\Bbb Z)\subset\Gamma_0$.
 In particular, we get that $Ae_{i,j}A^{-1}\in$ ${\Bbb M}_n(\Bbb Q)$, for all $i\not= j$, where $e_{i,j}$ denotes the elementary matrix $(e_{i,j})_{k,l}=\delta_{(i,j),(k,l)}$. Thus, if $A=(a_{i,j})_{1\leq i,j\leq n}$ and $A^{-1}=(b_{i,j})_{1\leq i,j\leq n}$, then we deduce that $a_{k,i}b_{j,l}\in \Bbb Q,$ for all $1\leq i,j,k,l\leq n$ with $i\not=j.$ Let $l_1,l_2,k_3$ such that $b_{j,l_j}\not=0$, for $j\in\{1,2\}$ and $a_{k_3,3}\not=0$. Thus, $a_{k,i}\in\Bbb Q b_{1,l_1}^{-1},$ for all $1\leq k,i\leq n$ with $i\not=1$ and $a_{k,i}\in \Bbb Qb_{2,l_2}^{-1}$, for all $1\leq k,i\leq n$ with $i\not=2$. Since we also have that $a_{k_3,3}b_{1,l_1},a_{k_3,3}b_{2,l_2}\in\Bbb Q$ and $a_{k_3,3}\not=0$, we deduce that $\Bbb Qb_{1,l_1}^{-1}=\Bbb Qb_{2,l_2}^{-1}$. Altogether, it follows that we can find $b\not=0$ such that $a_{k,i}\in \Bbb Qb^{-1}$, for all $1\leq k,i\leq n$, hence by replacing $A$ with $bA$ we can assume that $A=(a_{i,j})_{1\leq i,j\leq n}\in SL(n,\Bbb Q)$. Finally, let $s\in\Bbb N$ such that $sa_{i,j}\in\Bbb Z$, for all $1\leq i,j\leq n$ and replace $A$ with $sA$. \hfill$\blacksquare$
\vskip 0.1in
Recall that two ergodic actions $\Gamma\curvearrowright (X,\mu)$ and $\Lambda\curvearrowright (Y,\nu)$ are called {\bf stable orbit equivalent} if there exists an orbit equivalence $\theta:\Cal R_{\Gamma}\rightarrow\Cal R_{\Lambda}^t$, for some $t>0$. In the case $t=1$, we say that these actions are {\bf orbit equivalent}.
\vskip 0.1in
\proclaim {5.3. Corollary} Let $n\geq 3$ and let $\Gamma\subset$ SL$_n(\Bbb Z)$ be a finite index subgroup. Let $I_1,I_2$ be two infinite increasing sequences of prime numbers such that $|I_1\Delta I_2|=\infty$. Then the actions $\alpha_{I_1}$ and $\alpha_{I_2}$ are not stable orbit equivalent. In particular, if $\{I_t\}_{t\in \Bbb R}$ is a family of infinite increasing sequences of prime numbers such that $|I_s\Delta I_t|=\infty$, for all $s\not=t$, then $\{\alpha_{I_t}\}_{t\in \Bbb R}$ is an uncountable family of non stable orbit equivalent profinite actions of $\Gamma$. 
\endproclaim
{\it Proof.} Assume that  the  actions $\alpha_{I_1}$ and $\alpha_{I_2}$ are
 stable orbit equivalent. Note that $\alpha_{I_1}$ is an ergodic, profinite action
 of a property (T) group ([K]). Thus, Theorem 4.1. implies that we can find two 
finite index subgroups $\Gamma_1,\Gamma_2\subset \Gamma$ and ergodic components
 $Y_j\subset X_{I_j}$ for $\Gamma_j$, for all $j\in\{1,2\}$, such that the 
actions $\Gamma_1\curvearrowright Y_1$ and $\Gamma_2\curvearrowright Y_2$ are
 conjugate. Proposition 5.1. then gives that $|I_1\Delta I_2|<\infty$.
\hfill$\blacksquare$
\vskip 0.05in
\noindent
{\bf 5.4. Remarks.}
(1). Note that the first construction of an uncountable family of non stable orbit equivalent profinite actions for SL$_n(\Bbb Z)$ was obtained by S.L. Gefter and V.Y. Golodets ([GG]) as a consequence of Zimmer's cocycle superrigidity result ([Z]). 
   Later on, S. Thomas proved that the actions from [GG] moreover produce equivalence relations which are incomparable with respect to Borel reducibility ([T]).

(2). Recently, N. Ozawa and S. Popa showed that the Corollary 5.3. also holds true $n=2$, i.e.  every finite index subgroup $\Gamma$ of SL$_2(\Bbb Z)$ (e.g. $\Gamma=\Bbb F_m$, $2\leq m<\infty$) admits uncountably many non orbit equivalent profinite actions ([OP]). 

\vskip 0.1in

Next, we indicate how to construct non orbit equivalent profinite actions for the groups SL$_n(\Bbb Z)\ltimes\Bbb Z^n$ ($n\geq 2$). 
To do this, we first introduce a new isomorphism invariant for measurable equivalence relations. 
Note that in [P5], S. Popa studied a related version of this invariant for actions of countable groups on the hyperfinite II$_1$ factor.

\vskip 0.1in
\noindent
{\bf 5.5. Definition.} Let $\Cal R$ be a countable ergodic measure preserving
 equivalence relation. Then we let $\Cal S(\Cal R)$  be the set of $t>0$ such that $\Cal R^t$ is induced by a free ergodic measure preserving action of a countable group.
\vskip 0.1in
\noindent
{\bf 5.6. Remarks.}
(0). We gather in this remark a few elementary properties of the $\Cal S$-invariant. By definition, we have that $\Cal S(\Cal R^t)=\Cal S(\Cal R)/t$, for all $t>0$.
Moreover,  if $t\in \Cal S(\Cal R)$ and $n\geq 1$ is an integer, then $nt\in\Cal S(\Cal R)$. Indeed, let $\Cal R$ be an equivalence relation
induced by a free ergodic measure preserving action $\Gamma\curvearrowright (X,\mu)$. Then $\Cal R^n$ can be realized as the equivalence relation induced by the product action $\Gamma\times \Bbb Z_n\curvearrowright (X\times \Bbb Z_n,\mu\times c)$, where $c$ is the normalized counting probability measure on $\Bbb Z_n$ and $\Bbb Z_n$ acts on itself by group multiplication. 
Also, note that if $\Cal F(\Cal R)=\{t>0|\Cal R^t\simeq\Cal R\}$ denotes the fundamental group of $\Cal R$ and if $t\in\Cal S(\Cal R)$, then $t\Cal F(\Cal R)\subset\Cal S(\Cal R)$.
Finally, if $\Cal R$ is an ergodic hyperfinite equivalence relation, then $\Cal S(\Cal R)=\Bbb R_{+}^*$

(1). Feldman and Moore asked whether every  countable ergodic measure preserving  equivalence relation is induced by a free action of a countable group ([FM]). This  question was answered in the negative  by A. Furman, who  showed that if $n\geq 3$, then the equivalence relation $\Cal R$ induced by the action SL$_n(\Bbb Z)\curvearrowright\Bbb T^n$ satisfies $\Cal S(\Cal R)\subset \Bbb Q_{+}^*$ ([Fu2]).

(2). More examples of such equivalence relations  arise from the work of  S. Popa who showed that the equivalence relation $\Cal R$  induced by a Bernoulli action of a  property (T) group or of a product of two groups, one infinite and one non-amenable, satisfies $\Cal S(\Cal R)=\Bbb N^{*}$ ([P2,3]).
  Even more surprising, there are equivalence relations $\Cal R$ such that $\Cal S(\Cal R)=\emptyset$ ([Fu2],[P2]).
\vskip 0.05in

The following result is a direct consequence of Theorem 4.1.
\proclaim{5.7. Corollary} Let $\Gamma\curvearrowright^{\alpha}X$ be as in Theorem A and let $\Cal R_{\Gamma}$ be the equivalence relation induced by $\alpha$.  Suppose that $\alpha$ is the limit of the actions $\Gamma\curvearrowright^{\alpha_n} X_n$, with $X_n$ finite. Then $\Cal S(\Cal R_{\Gamma})=\{m/|X_n||m,n\geq 1\}\subset \Bbb Q_{+}^*.$
\endproclaim
{\it Proof.} Let $t\in\Cal S(\Cal R_{\Gamma})$. Thus there exists a free ergodic measure preserving action of a countable group $\Lambda$ such that $\Cal R_{\Gamma}^t=\Cal R_{\Lambda}$. By applying Theorem 4.1., we get that $t=[\Lambda:\Lambda_0]/|X_n|$, for some $n$ and some finite index subgroup $\Lambda_0$ of $\Lambda$. This shows that $\Cal S(\Cal R_{\Gamma})\subset\{m/|X_n||m,n\geq 1\}.$ For the other inclusion, note that by Remark 5.6.(0). we only need to show that $1/|X_n|\in\Cal S(\Cal R_{\Gamma})$, for all $n$.  This is clear, since for any $a\in X_n$, the restriction of $\Cal R_{\Gamma}$ to $X_{a,n}$ (which has measure $1/|X_n|$) is induced by the action $\Gamma_{a,n}\curvearrowright X_{a,n}$. 
\hfill$\blacksquare$
\vskip 0.05in

It is natural to ask to characterize the subsets of $\Bbb R_{+}^*$ which can appear as the $\Cal S$-invariant of a measurable equivalence relation $\Cal R$.
Next, we specialize  Corollary 5.7. to some concrete actions of SL$_n(\Bbb Z)\ltimes\Bbb Z^n$ ($n\geq 2$) and  deduce that the invariant $\Cal S(\Cal R)$ can equal $\Bbb Q_{+}^*$ as well as many interesting arithmetic sets. 

\vskip 0.1in
Fix $n\geq 2$. Let $J=\{d_i\}_{i\geq 1}$ be a set of integers such that $d_i<d_{i+1}$ and $d_i|d_{i+1}$ (equivalently, such that $d_{i+1}\Bbb Z\subsetneq d_i\Bbb Z$), for all $i\geq 1$.
 Notice that, for every $i\geq 1,$ the subgroup $(d_i\Bbb Z)^n$ of $\Bbb Z^n$ is left invariant	 by the action of SL$_n(\Bbb Z)$. Thus, following 1.6. we can construct an ergodic profinite action SL$_n(\Bbb Z)\ltimes\Bbb Z^n\curvearrowright^{\beta_J}\varprojlim (\Bbb Z/d_i\Bbb Z)^n$. 
This action is also free. Indeed, we have that $\cap_{i}(d_i\Bbb Z)^n=\{0\}$. Also, if $A\in$ SL$_n(\Bbb Z)\setminus\{I\}$, then $\{x\in \Bbb Z^n|Ax=x\}$ is an infinite index subgroup of $\Bbb Z^n$ and thus Lemma 1.7. implies that $\beta_J$ is free.

Moreover, for any set $J$ as above, the action $\beta_J$ satisfies the assumptions of Theorem 5.7. This is true since for all $n\geq 2$, the pair (SL$_n(\Bbb Z)\ltimes\Bbb Z^n,\Bbb Z^n)$ has relative property (T) ([K]), SL$_n(\Bbb Z)$ is finitely generated and the restriction ${\beta_J}_{|\Bbb Z^n}$ is ergodic.

\proclaim {5.8. Corollary}
Let $J$ be a set of integers as above and denote by $\Cal R_{J}$ the equivalence relation 
induced by $\beta_J$. Then $\Cal S(\Cal R_{J})=\{k/d_i^n|k,i\geq 1\}$. In particular,
 $\Cal S(\Cal R_{J})=\Bbb Q^{*}_{+}$ for some  set $J$.  Also, 
if $J_1=\{d_{i}^1\}_{i\geq 1}$ and $J_2=\{d_{i}^2\}_{i\geq 1}$ are two sets as 
above such that there exists $i$ with either $d_{i}^2\not| d_{j}^1$, for all $j$, 
or $d_{i}^1\not| d_{j}^2$,  for all $j$, then the profinite actions $\beta_{J_1}$ 
and $\beta_{J_2}$ are not orbit equivalent.  
\endproclaim
Note that Remark 5.5. and Corollary 5.6. provide several examples of actions $\Gamma\curvearrowright X$ of property (T) groups $\Gamma$ such that $\Cal S(\Cal R_{\Gamma})$ is countable. We prove below that this is in fact true for any property (T) group $\Gamma$, regardless of the action $\Gamma\curvearrowright X$. 
The proof we give is obtained by assembling together ideas from Section 4 in [P4]. We refer
 the reader to Section 5 in [P6] for the definition of property (T) for 
equivalence relations.  

\proclaim{5.9. Theorem} Let  $\Cal R$ be a countable ergodic measure preserving equivalence relation on a standard probability space $(X,\mu)$.
 If $\Cal R$ has property (T), then $\Cal S(\Cal R)$ is countable.
\endproclaim
{\it Proof.} If $t\in\Cal S(\Cal R)$, let $\Gamma$ be a countable group such that $\Cal R^t$ is induced by a free action of $\Gamma$. Since $\Cal R$ has property (T), $\Cal R^t$ has property (T), hence $\Gamma$ has property (T) ([Fu1]).
Also, recall that by a theorem of Y. Shalom ([Sh2]), any property (T) group is the quotient of a finitely presented, property (T) group.

Assume by contradiction that $\Cal S(\Cal R)$ is uncountable.  Using the preceding discussion we deduce that there exists a property (T) group $\Gamma$ such that the set $\Cal S$ of $t>0$ for which $\Cal R^t$ is induced by a free action of a quotient of $\Gamma$ is uncountable. 
Since $\Cal S$ is uncountable, after eventually replacing $\Cal R$ with an amplification of it, we can assume that $\Cal S\cap (3/4,1)$ is uncountable. 
Let $\{X_t\}_{t\in [0,1]}$ be measurable sets such that $X_t\subset X_{t'}\subset X$ and $\mu(X_t)=t$ for all $0\leq t\leq t'\leq 1$. For every $t\in S$ let $\alpha_t$ be the action of $\Gamma$ on $X_t$ (through one of its quotients) such that $$\Cal R\cap (X_t\times X_t)=\{(x,\alpha_t(\gamma)(x)|x\in X_t,\gamma\in\Gamma\}.$$ Extend $\alpha_t$ to $X$ by letting $\Gamma$ act identically on $X\setminus X_t$. Next, endow $\Cal R$ with the measure $\tilde\mu(C)=\int_{X}|\{y|(x,y)\in C\}|d\mu(x)$, for any Borel subset $C$ of $\Cal R$ ([FM]).
For every $t,t'\in S\cap (0,1]$ let $\pi_{t,t'}:\Gamma\rightarrow\Cal U(L^2(\Cal R,\tilde\mu))$ be the unitary representation defined by $$\pi_{t,t'}(\gamma)(f)(x,y)=f(\alpha_t(\gamma^{-1})(x),\alpha_{t'}(\gamma^{-1})(y)),$$ for all $\gamma\in\Gamma,f\in L^2(\Cal R,\tilde\mu),(x,y)\in\Cal R$.

Using the property (T) of $\Gamma$  as in Section 2 of [Hj], we can find $t,t'\in S\cap (3/4,1)$ such that $t<t'$ and $\pi_{t,t'}$ admits an invariant vector $f\in L^2(\Cal R,\tilde \mu)$ with $||f-1_{\Delta}||_{L^2(\Cal R,\tilde\mu)}< 1/4$, where $\Delta=\{(x,x)|x\in X\}$. Set $$A:=\{x\in X|\exists ! y=\phi(x)\in X, |f(x,y)|>1/2\},$$ then,  since $f$ is $\pi_{t,t'}$-invariant, $A$ is $\alpha_t(\Gamma)$-invariant. If $$A_0:= \{x\in X||f(x,x)-1|^2+\sum_{y\Cal Rx,y\not=x}|f(x,y)|^2<1/4\},$$ then $A_0\subset A$ and $\mu(X\setminus A_0)/4\leq ||f-1_{\Delta}||_{L^2(\Cal R,\tilde\mu)}^2<1/16.$ Altogether, we get that $\mu(A)>3/4$. Since the restriction of $\alpha_t$ to $X_t$ is ergodic and $\mu(X_t)>1/2$ we deduce that $X_t\subset A$.

Now, let $B=\phi(X_t)=\{y\in X|\exists x\in X_t, |f(x,y)|>1/2\}$, then $\mu(B)\leq \mu(X_t).$ 
Proceeding as above, if we denote $B_0:= \{y\in X||f(y,y)-1|^2+\sum_{y\Cal Rx,y\not=x}|f(x,y)|^2\leq 1/4\},$ then $\mu(B_0)>3/4$
  and since $B_0\cap X_t\subset B$, we deduce that $\mu(B)>\mu(X_t)-1/4>1/2.$
 Moreover, since $B$ is $\alpha_{t'}(\Gamma)$-invariant and since the restriction of $\alpha_{t'}$ to $X_{t'}$ is ergodic, we get  that $X_{t'}\subset B$.
Thus, $\mu(B)\geq\mu(X_{t'})$ a contradiction to $\mu(B)\leq\mu(X_t)$.\hfill$\blacksquare$
\vskip 0.1in

\head \S 6. Compact and weakly compact actions.\endhead
\vskip 0.1in
The aim of this Section is to investigate N. Ozawa and S. Popa's recent notion of weak compactness for actions ([OP]) and its relation with the classical notion of compactness for actions.
  Thus, we show  that any compact action is  weakly compact, that weak compactness is an orbit equivalence invariant property and that weak compactness passes to quotients. Note that these facts have been obtained independently in [OP] (see Propositions 3.2. and 3.4.)
 As a corollary, it follows that if  $\sigma$ is a compact action (e.g. profinite) and $\alpha$ is an action which is orbit equivalent to $\sigma$, then any quotient $\beta$ of $\alpha$ does not have stable spectral gap, in the sense on [P3]. This extends the main result of [ET] in the case of orbit equivalence.  We end the Section by giving a new characterization of  compact actions, in terms of the von Neumann algebras associated (Theorem 6.9.). 
\vskip 0.05in 
\noindent
{\bf 6.1. Compact actions.} An ergodic measure preserving action $\Gamma\curvearrowright^{\sigma}(X,\mu)$ is called {\bf compact} (or {\bf isometric}) if one of the following equivalent conditions holds true (see [G] for a reference):
\vskip 0.04in

$(i)$ The closure  of $\Gamma$ in Aut$(X,\mu)$ is compact.
\vskip 0.03in
$(ii)$ The induced unitary representation $\pi:\Gamma\rightarrow\Cal U(L^2(X,\mu))$ decomposes as a direct sum of finite dimensional representations.
\vskip 0.03in
$(iii)$ $\sigma$ is of the form $\Gamma\curvearrowright (K/L,\tilde m)$ and  is defined by $\gamma (kL)=\alpha(\gamma)kL$, where $K$ is a compact group, $L\subset K$ is a closed subgroup, $\tilde m$ is the projection of the Haar measure $m$ of $K$ onto $K/L$ and $\alpha:\Gamma\rightarrow K$ is a homomorphism with a dense image.
\vskip 0.05in
Note that every ergodic profinite action is compact. Indeed, if an action $\Gamma\curvearrowright X$ is the limit of actions $\Gamma\curvearrowright X_n$, with $X_n$ finite, then the unitary representation $\Gamma\curvearrowright L^2X$ is the direct sum of the finite dimensional representations $\Gamma\curvearrowright L^2X_n\ominus L^2X_{n-1}$. 
\vskip 0.1in
Next, we introduce some new notation, that we will maintain throughout this section.
If $\Gamma\curvearrowright^{\sigma} (X,\mu)$ is a measure preserving action, then we denote by  $\Gamma\curvearrowright^{\tilde\sigma}(X\times X,\mu\times\mu)$ the double action $\tilde\sigma=\sigma\times\sigma$.  Also, we denote by $\tau:L^1(X,\mu)\rightarrow \Bbb C$ the integral with respect to $\mu$. Then $\tau\otimes \text{id}:L^1(X\times X,\mu\times\mu)\rightarrow L^1(X,\mu)$ (respectively $\text{id}\otimes\tau:L^1(X\times X,\mu\times\mu)\rightarrow L^1(X,\mu)$) are the conditional expectations onto $1\otimes L^1(X,\mu)$ (respectively onto $L^1(X,\mu)\otimes 1$). Finally, for every measurable set $A\subset X$, we denote $\tilde A=A\times (X\setminus A)\subset X\times X$.
\vskip 0.1in
\noindent
{\bf 6.2. Weakly compact actions}.  An ergodic measure preserving action $\Gamma\curvearrowright^{\sigma}(X,\mu)$ is called {\bf weakly compact} ([OP]) if there exist a sequence $\{\eta_n\}_{n\geq 1}\in L^1(X\times X,\mu\times\mu)$ of positive vectors with $||\eta_n||_1=1$  such that 
\vskip 0.04in
$(i)$  $\lim_{n\rightarrow\infty}||1_{\tilde A}\eta_n||_1=0,$ for every measurable set $A\subset X$.
\vskip 0.04in
$(ii)$  $\lim_{n\rightarrow\infty}||\eta_n-\eta_n\circ{\tilde\sigma}(\gamma)||_1=0,$ for all $\gamma\in\Gamma$.
\vskip 0.04in
$(iii)$  $(\tau\otimes\text{id})(\eta_n)=(\text{id}\otimes\tau)(\eta_n)=1,$ for all $n\geq 1.$

\vskip 0.1in

The class of weakly compact actions has been introduced by N. Ozawa and S. Popa in the course to establishing the following remarkable result: if $\Bbb F_m\curvearrowright^{\sigma}(X,\mu)$ is a free ergodic profinite action of a free group $\Bbb F_m$, $2\leq m\leq\infty$, then $L^{\infty}(X,\mu)$ is the unique Cartan subalgebra of the crossed product von Neumann algebra $L^{\infty}(X,\mu)\rtimes_{\sigma}\Bbb F_m$, up to unitary conjugacy ([OP]).

\proclaim {6.3. Proposition} An ergodic measure preserving action $\Gamma\curvearrowright^{\sigma} (X,\mu)$ is compact {\it iff}
there exist a sequence $\{\eta_n\}_{n\geq 1}\in L^1(X\times X,\mu\times\mu)$ of $\tilde\sigma(\Gamma)$-invariant, positive vectors satisfying $||\eta_n||_1=1$, for all $n$, and  
  $\lim_{n\rightarrow\infty}||1_{\tilde A}\eta_n||_1=0,$ for every measurable set $A\subset X$. In particular, any compact action is weakly compact.
\endproclaim

{\it Proof.} ($\Longrightarrow$) Assume first that $\sigma$ is compact. Thus, using 6.1.$(ii)$, we can identify
 $(X,\mu)=(K/L,\tilde m)$ where $\Gamma$ acts on $K/L$ via a homomorphism $\alpha:\Gamma\rightarrow K$. Let $d$ be a metric on $K/L$ which is invariant under the left $K$-action.
For every $n$, define $A_{n}=\{(x,y)\in X\times X|d(x,y)<1/n\}.$ Then $(\mu\times\mu)(A_{n})>0$ and it is routine to check that the vectors $\eta_{n}=1_{A_{n}}/(\mu\times \mu)(A_{n})$ verify the desired conditions.

\vskip 0.05in

($\Longleftarrow$) Conversely, assume that there exist vectors $\eta_n\in L^1(X\times X,\mu\times\mu)$ satisfying the hypothesis and suppose by contradiction that $\sigma$ is not a compact action.  Let $\Gamma\curvearrowright^{\alpha} (Y,\mu)$ be the maximal compact quotient of $\sigma$ with the quotient map $p:(X,\mu)\rightarrow (Y,\nu)$. 
If we denote by $\pi:\Gamma\rightarrow\Cal U(L^2(X,\mu))$ the unitary representation induced by $\sigma$, then $L^2(Y,\nu)$ (viewed as a Hilbert subspace of $L^2(X,\mu)$, via $p$) is precisely the closure of the union of all finite dimensional, $\pi(\Gamma)$-invariant Hilbert subspaces of $L^2(X,\mu)$. This implies (by using the standard identification of $L^2(X\times X,\mu\times\mu)$ with the Hilbert-Schmidt operators on $L^2(X,\mu)$)  that if $\eta\in L^2(X\times X,\mu\times\mu)$ is a $\tilde\sigma(\Gamma)$-invariant vector, then $\eta\in L^2(Y\times Y,\nu\times\nu)$ ([G]). In particular, we derive that $\eta_n\in L^1(Y\times Y,\nu\times\nu)$, for all $n$.

Next, denote by $E:L^{\infty}(X,\mu)\rightarrow L^{\infty}(Y,\nu)$ and by $\tilde E:L^{\infty}(X\times X,\mu\times\mu)\rightarrow L^{\infty}(Y\times Y,\nu\times\nu)$ the conditional expectations onto $L^{\infty}(Y,\nu)$ and onto $L^{\infty}(Y\times Y,\nu\times\nu)$, respectively. We claim that there exists a measurable set $A\subset X$ such that $\tilde E(1_{\tilde A})\geq (1/6) 1_{Y\times Y}.$ 
Indeed, since $\sigma$ is ergodic, by Rokhlin's skew product theorem ([G])  we can decompose $(X,\mu)=(Y,\nu)\times (Z,\rho)$, where either $Z=[0,1]$ and $\rho$ is the Lebesgue measure or $Z=\{1,2,..,n\}$ and $\rho=1/n\sum_{i=1}^n\delta_{i}$, for some $n\geq 2$ (here we use the fact that $X\not=Y$). 
In the first case, if we let $A=Y\times [0,1/2]$, then $E(1_A)=E(1_{X\setminus A})=(1/2) 1_Y$.
In the second case, if we let $A=Y\times \{1,..,[n/2]\}$, then $E(1_A)=([n/2]/n) 1_Y\geq (1/3)1_Y$ and $E(1_{X\setminus A})=(1-[n/2]/n) 1_Y\geq (1/2)1_Y$. Finally, since $\tilde E(1_{\tilde A})=E(1_A)E(1_{X\setminus A})$, we obtain the claim in both cases.

To end the proof,  just note that the above imply that $$||1_{\tilde A}\eta_n||_1=\int_{X\times X}1_{\tilde A}\eta_n \text{d}(\mu\times\mu)=\int_{Y\times Y}\tilde{E}(1_{\tilde A})\eta_n \text{d}(\nu\times\nu)\geq$$ $$(1/6)\int_{Y\times Y}\eta_n \text{d}(\nu\times\nu)=(1/6)||\eta_n||_1=1/6,\forall n,$$in contradiction with our assumption.\hfill$\blacksquare$

\vskip 0.05in
\noindent
{\bf 6.4. Remark}.  We next show (Theorem 6.6.) that the class of weakly compact actions is closed under orbit equivalence (see also [OP]). Note that this is not true for the class of compact actions.  To see this, recall that by Dye's theorem ([Dy]) the profinite (hence compact) action $\Bbb Z\curvearrowright \varprojlim \Bbb Z/2^n\Bbb Z$ is orbit equivalent to the weakly mixing (hence not compact) Bernoulli action $\Bbb Z\curvearrowright [0,1]^{\Bbb Z}$. 

However, if we restrict to compact actions $\Gamma\curvearrowright^{\sigma}(X,\mu)$ of property (T) groups $\Gamma$, then this class is closed under orbit equivalence.
Indeed, if $\Lambda\curvearrowright^{\alpha}(Y,\nu)$ is a free ergodic measure preserving action which is orbit equivalent to $\sigma$, then by [Fu1], $\Lambda$  has property (T), while by Theorem 6.6. below, $\alpha$ is weakly compact.
Let $\eta_n\in L^1(Y\times Y,\nu\times\nu)$ be vectors satisfying Definition 6.2. Then $\eta_n^{1/2}\in L^2(Y\times Y,\nu\times\nu)$ forms a sequence of almost invariant vectors for the unitary representation induced by $\tilde\alpha$. Since $\Lambda$ has property (T), we can find a sequence of $\tilde\alpha$-invariant unit vectors  $\xi_n$ such that $\lim_{n\rightarrow\infty}||\eta_n^{1/2}-\xi_n||_2=0$.
Thus, if we denote $\eta_n'=\xi_n^2$, for all $n$, then $\eta_n'$ is a $\tilde\alpha$-invariant positive vector with $||\eta_n'||_1=1$ and $\lim_{n\rightarrow\infty}||\eta_n-\eta_n'||_1=0$. In particular, since $\lim_{n\rightarrow\infty}||1_{\tilde A}\eta_n||_1=0$, we get that $\lim_{n\rightarrow\infty}||1_{\tilde A}\eta_n'||_1=0$, for every measurable set $A\subset X$. Proposition 6.3. then implies that $\alpha$ is an isometric action.   

Note that A. Furman proved that in fact compact actions of property (T) groups are OE superrigid , i.e., in the above context, $\alpha$ must be (virtually) conjugate to $\sigma$ (personal communication).
 
\proclaim {6.5. Lemma} Let $\Gamma\curvearrowright^{\sigma}(X,\mu)$ be a weakly compact action and let $\eta_n\in L^1(X\times X,\mu\times\mu)$ be a sequence of vectors satisfying the conditions of Definition 6.2. Let $\phi\in$ Aut$(X,\mu)$ such that $\phi(x)\in\Gamma x$, a.e. $x\in X$ (i.e. $\phi\in [\Cal R_{\Gamma}]$), and denote $\tilde \phi=\phi\times\phi\in$ Aut$(X\times X,\mu\times \mu)$. Then $\lim_{n\rightarrow\infty}||\eta_n-\eta_n\circ\tilde\phi||_1=0.$
\endproclaim
{\it Proof.} Fix an enumeration $\{\gamma_1,\gamma_2,...\}$ of $\Gamma$. For every $i,k\geq 1$, let $A_i=\{x\in X|\phi(x)={\gamma_i}^{-1}x\}$, $S_k=\cup_{i=1}^k(A_i\times A_i)$ and $X_k=X\setminus(\cup_{i=1}^k A_i)$. Then for all $n$ and $k$ we have that  $$||\eta_n-1_{S_k}\eta_n||_1\leq\sum_{1\leq i\not= j\leq k}||1_{A_i\times A_j}\eta_n||_1+||1_{X\times X_k}\eta_n||_1+||1_{X_k\times X}\eta_n||_1\tag 6.5.a.$$
Now, condition 6.2.$(iii)$ implies that $||1_{X\times X_k}\eta_n||_1=||1_{X_k\times X}\eta_n||_1=\mu(X_k),$ for all $k$. Also, for all $i\not= j$, we have that $||1_{A_i\times A_j}\eta_n||_1\leq ||1_{\tilde{A_i}}\eta_n||_1$, thus by condition 6.2.$(i)$, we get that $\lim_{n\rightarrow\infty}||1_{A_i\times A_j}\eta_n||_1=0.$ By combining these facts with (6.5.a.) we get that $$\limsup_{n\rightarrow\infty}||\eta_n-1_{S_k}\eta_n||_1\leq 2\mu(X_k),\forall k\tag 6.5.b.$$

Next, for all $n$ and $k$, by triangle's inequality we get that $$||\eta_n-\eta_n\circ\tilde\phi||_1\leq 2||\eta_n-1_{S_k}\eta_n||_1+||1_{S_k}\eta_n-1_{S_k}\eta_n\circ\tilde\phi||_1\leq \tag 6.5.c.$$ $$2||\eta_n-1_{S_k}\eta_n||_1+\sum_{i=1}^k||\eta_n-\eta_n\circ\tilde\sigma(\gamma_i)||_1.$$
Since by condition 6.2.($ii$). we have that  $\lim_{n\rightarrow\infty}||\eta_n-\eta_n\circ\tilde\sigma(\gamma_i)||_1=0$, for all $i$, by using (6.5.b.) and (6.5.c.) together we derive that $\limsup_{n\rightarrow\infty}||\eta_n-\eta_n\circ\tilde\phi||_1\leq 4\mu(X_k),$ for all k. Finally, since $\mu(X_k)\rightarrow 0$, as $k\rightarrow\infty$, we get the conclusion. \hfill$\blacksquare$

\proclaim {6.6. Theorem} Let $\Gamma\curvearrowright^{\sigma}(X,\mu)$ be a weakly compact action and let $\Lambda\curvearrowright^{\alpha}(Y,\nu)$, $\Gamma\curvearrowright^{\beta}(Z,\rho)$ be two measure preserving actions.
\vskip 0.05in
$(i)$ If $\alpha$ is orbit equivalent to $\sigma$, then $\alpha$ is weakly compact.
\vskip 0.05in
$(ii)$ If $\beta$ is a quotient of $\sigma$, then $\beta$ is weakly compact.
\endproclaim
{\it Proof.} $(i)$. This is an immediate consequence of Lemma 6.5.

$(ii)$. Let $\eta_n\in L^1(X\times X,\mu\times \mu)$ satisfying Definition 6.2. Let $p:(X,\mu)\rightarrow (Z,\rho)$ be a $\Gamma$-equivariant quotient map and let $E:L^1(X\times X,\mu\times\mu)\rightarrow L^1(Z\times Z,\rho\times\rho)$ be the conditional expectation. Then $E$ is integral preserving, thus $\xi_n:=E(\eta_n)\in L^1(Y\times Y,\nu\times \nu)$ is a positive vector with $||\xi_n||_1=1$, for all $n$. Also, we have that $||E(\eta)||_1\leq ||\eta||_1$ and  $E(\gamma \eta)=\gamma E(\eta)$, for all $\eta\in L^1(X\times X,\mu\times\mu)$ and every $\gamma\in\Gamma$. Thus $$||\xi_n-\xi_n\circ{\tilde{\beta}}(\gamma)||_1=||E(\eta_n-\eta_n\circ\tilde{\sigma}(\gamma))||_1\leq ||\eta_n-\eta_n\circ\tilde{\sigma}(\gamma)||_1,$$ hence the $\xi_n$'s verify condition 6.2.$(ii)$.
Moreover, if $A\subset Y$ is a measurable set, then $$||1_{\tilde A}\xi_n||_1=||E(1_{\widetilde{p^{-1}(A)}}\eta_n)||_1\leq ||1_{\widetilde{p^{-1}(A)}}\eta_n||_1,$$ thus the $\xi_n$'s also verify condition 6.2.$(i)$. Since condition 6.2.$(iii)$. is clearly verified, we deduce that $\beta$ 
is  weakly compact.\hfill$\blacksquare$
\vskip 0.1in

 Let $\Gamma\curvearrowright^{\sigma}(X,\mu)$ be a measure preserving action and denote by $\pi$ the induced unitary representation of $\Gamma$ on  $L^2(X,\mu)$ and by $\pi_0$ its restriction to $L^2(X,\mu)\ominus\Bbb C1$. Then we say that $\sigma$ has {\bf stable spectral gap} if the representation $\pi_0\otimes\pi_0$ does not weakly contain the trivial representation (see Section 3 in [P3]). 
\vskip 0.05in
\noindent
\proclaim {6.7. Proposition}  Any weakly compact action $\Gamma\curvearrowright^{\sigma} (X,\mu)$ on a standard probability space $X$ does not have  stable spectral gap.
\endproclaim
\vskip 0.05in
{\it Proof.} Let $\eta_n\in L^1(X\times X,\mu\times\mu)$ be a sequence of vectors satisfying Definition 6.2. Then $\xi_n=\eta_n^{1/2}\in L^2(X\times X,\mu\times\mu)=L^2(X,\mu)\overline{\otimes}L^2(X,\mu)$ form a sequence of almost invariant  unit vectors for $\pi\otimes\pi$. Also $\xi_n$ verify $$\lim_{n\rightarrow\infty}||1_{\tilde A}\xi_n||_2=0\tag 6.7.a.$$ for any measurable set $A\subset X$.
For all $n$, denote by $\xi_n^0$ (resp. $\xi_n^1$ and $\xi_n^2$) the orthogonal projection of $\xi_n$  onto the Hilbert space $(L^2(X,\mu)\ominus\Bbb C1)\overline{\otimes}(L^2(X,\mu)\ominus\Bbb C1)$ (resp. onto $L^2(X,\mu)\otimes 1$ and onto $1\otimes (L^2(X,\mu)\ominus\Bbb C1)$). Then $\xi_n^0$ form a sequence of almost invariant vectors for $\pi_0\otimes\pi_0$. Thus, to prove the proposition it suffices to show that $||\xi_n^0||_2\nrightarrow 0$ as $n\rightarrow\infty$. 

Now, fix $k\geq 1$ and let $\{A_1,A_2,..,A_k\}$ be a measurable partition of $X$ such that $\mu(A_i)=1/k$, for all $i\in\{1,..,k\}$. Then (6.7.a.) easily implies that $\lim_{n\rightarrow\infty}||\xi_n-(\sum_{i=1}^k1_{A_i\times A_i})\xi_n||_2=0,$ thus $$\lim_{n\rightarrow\infty}||(\sum_{i=1}^k1_{A_i\times A_i})\xi_n||_2=1\tag 6.7.b.$$
Also, it is easy to check that if $j\in\{1,2\}$, then $$||(\sum_{i=1}^k1_{A_i\times A_i})\xi_n^j||_2=\sqrt{1/k}||\xi_n^j||_2\tag 6.7.c.$$ Since $\xi_n=\xi_n^0+\xi_n^1+\xi_n^2$, for all $n$, by combining (6.7.b.) and (6.7.c.) we get that $\limsup_{n\rightarrow\infty}||\xi_n^0||_2\geq\limsup_{n\rightarrow\infty} ||(\sum_{i=1}^k1_{A_i\times A_i})\xi_n^0||_2\geq 1-2\sqrt{1/k}.$ In particular, for $k=5$, we get that  $\limsup_{n\rightarrow\infty}||\xi_n^0||_2\geq 1-2/{\sqrt{5}}>0$. \hfill$\blacksquare$

\vskip 0.1in
To give some context for the next result,  let $\Lambda\curvearrowright^{\alpha}Y=Z\times V$  be a skew product action, i.e. an action given by $\lambda (z,v)=(\lambda z,w(\lambda,z)v)$, for an action $\Lambda\curvearrowright^{\beta} Z$ and a cocycle $w:\Lambda\times Z\rightarrow$ Aut$(V)$. Then the main result of [ET] says that if
$\Gamma\curvearrowright^{\sigma} X$ is a profinite action  such that there exists a countable-to-one Borel map $\theta:X\rightarrow Y$ with $\theta(\Cal R_{\Gamma})\subset\Cal R_{\Lambda}$ and if the image of $w$ in Aut$(V)$ is countable, then $\beta$ does not have stable spectral gap.
We note below that if we assume that $\theta$ is in fact an orbit equivalence between $\Cal R_{\Gamma}$ and $\Cal R_{\Lambda}$, then we ca generalize the above result ([ET]), by removing the countability assumption on $w$.

\proclaim {6.8. Corollary} Let  $\Gamma\curvearrowright^{\sigma}(X,\mu)$ be an ergodic compact action. Assume that $\Lambda\curvearrowright^{\alpha}(Y,\nu)$ is an action which is orbit equivalent to $\sigma$ and that $\Lambda\curvearrowright^{\beta}(Z,\rho)$ is a quotient of $\alpha$. Then $\beta$ does not have stable spectral gap.
\endproclaim
{\it Proof.} Just combine 6.6. and 6.7.
\vskip 0.1in
\noindent

We conclude this Section by giving a characterization of compact  actions in the language of von Neumann algebras. Let us first recall a few notions (see [P6] for all this). Let $(N,\tau)$ be a separable finite von Neumann algebra together with a normal, faithful trace $\tau$. Endow $N$ with the Hilbert norm $||x||_2=\tau(x^*x)^{1/2}$, for all $x\in N$, and let $L^2N$ be the completion of $N$ with respect to this norm. Let $B\subset N$ be a von Neumann subalgebra.
The {\bf quasi-normalizer} of $B$ in $N$, denoted $q\Cal N_{N}(B)$, is the set of $x\in N$ for which there exist $x_1,x_2,..,x_n\in N$ satisfying $xB\subset \sum Bx_i$ and $Bx\subset \sum x_iB$ (see 1.4.2. in [P6]). Then  $ q\Cal N_{N}(B)$ forms a $*$-algebra and the von Neumann algebra it generates, $q\Cal N_{N}(B)''$,  is a von Neumann subalgebra of $N$ which contains $B$. If $q\Cal N_{N}(B)''=N$, then we say that $B$ is {\bf quasi-regular} in $N$ ([P6]).

Also, we denote by $E_B:N\rightarrow B$ the conditional expectation onto $B$ and by $e_B:L^2N\rightarrow L^2B$ the orthogonal projection onto $L^2B$. Then {\bf Jones' basic construction} for the inclusion $B\subset N$ is denoted by $<N,e_B>$ and is defined as the von Neumann algebra generated by $N$ and $e_B$ inside $\Bbb B(L^2N)$. Note that $<N,e_B>$ is endowed with a natural semifinite trace given by Tr$(xe_By)=\tau(xy)$, for all $x,y\in N$. 

Next, recall that $N$ has {\bf property (H) (Haagerup) relative to } $B$ (see definition 2.1. in [P6]) if there exists a sequence of normal, $B$-bimodular, completely positive maps $\Phi_n:N\rightarrow N$ such that $\tau\circ\Phi_n\leq\tau$,  $\lim_{n\rightarrow\infty}||\Phi_n(z)-z||_2=0$, for all $z\in N$, and $\Phi_n$ is ''compact relative to $B$'': for every $\varepsilon>0$, there exists a projection $P\in <N,e_B>$ of finite trace such that the bounded operator $T_{\Phi_n}\in \Bbb B(L^2N)$ given by $T_{\Phi_n}(x)=\Phi_n(x)$, for all $x\in N$, satisfies $||T_{\Phi_n}(1-P)||\leq\varepsilon$. For example, in the case $N$ decomposes as a crossed product $B\rtimes_{\sigma}\Gamma$, for an action $\sigma:\Gamma\rightarrow$ Aut$(B,\tau)$ of a countable group $\Gamma$, we have that $N$ has property (H) relative to $B$ iff $\Gamma$ has Haagerup's property (see Proposition 3.1. in [P6]). Also, note that in the case $B=N$ relative property (H) amounts to property (H) of the single algebra $N$ (see 2.0.2. in [P6]).  

Finally, recall that to every measure preserving action $\Gamma\curvearrowright^{\sigma}(X,\mu)$ one associates the crossed product von Neumann algebra $L^{\infty}X\rtimes_{\sigma}\Gamma$ ([MvN]). Roughly speaking, this  algebra is generated by a copy of $L^{\infty}(X,\mu)$ and a copy of $\Gamma=\{u_{\gamma}\}_{\gamma\in\Gamma}$  subject to the relations $u_{\gamma}fu_{\gamma}^*=f\circ\sigma({\gamma}^{-1})$, for all $\gamma\in\Gamma$ and $f\in L^{\infty}X$.   Note that $L^{\infty}X\rtimes_{\sigma}\Gamma$ is a finite von Neumann algebra with the trace $\tau$ given by $\tau(\sum_{\gamma}f_{\gamma}u_{\gamma})=\int_{X}f_e d\mu$ and that the von Neumann subalgebra generated by $\{u_{\gamma}|\gamma\in\Gamma\}$ is isomorphic to the group von Neumann algebra $L\Gamma$.

\vskip 0.05in
\proclaim{6.9. Theorem} let $\Gamma\curvearrowright^{\sigma}(X,\mu)$ be an ergodic measure preserving  action. Then the following are equivalent
\vskip 0.03in
$(a)$ $\sigma$ is compact.
\vskip 0.03in
$(b)$ $L^{\infty}X\rtimes_{\sigma}\Gamma$ has property H relative to $L\Gamma$.
\vskip 0.03in
$(c)$ $L\Gamma$ is quasi-regular in $L^{\infty}X\rtimes_{\sigma}\Gamma$.
\endproclaim

{\it Proof.} $(a)\Longrightarrow (b)$ Assume that $\sigma$ is compact. Then, as in the proof of 6.3., we identify $\sigma$ with the action $\Gamma\curvearrowright K/L$, we let $d$ be a $K$-invariant metric on $K/L$ and we set $A_n=\{(x,y)\in X\times X|d(x,y)<1/n\}$, $\eta_n=1_{A_n}/(\mu\times\mu)(A_n)$, for all $n\geq 1$.  For every $n$, let $\phi_n:L^{\infty}X\rightarrow L^{\infty}X$ be given by $$\phi_n(f)(x)=\int_{X}\eta_n(x,y)f(y)d\mu(y), \forall f\in L^{\infty}X,x \in X.$$ Then  $\phi_n$ is a completely positive (since $\phi_n$ is positive and $L^{\infty}X$ is abelian), unital, integral preserving map which commutes with $\sigma$. 
 Thus, $\phi_n$ extends to a $L\Gamma$-bimodular unital completely positive map $\Phi_n:L^{\infty}X\rtimes_{\sigma}\Gamma\rightarrow L^{\infty}X\rtimes_{\sigma}\Gamma$ through the formula $\Phi_n(\sum_{\gamma}f_{\gamma}u_{\gamma})=\sum_{\gamma}\phi_n(f_{\gamma})u_{\gamma}$, where $f_{\gamma}\in L^{\infty}X$, for all $\gamma\in\Gamma$.  

Next, we claim that $\lim_{n\rightarrow\infty}||\Phi_n(z)-z||_2=0$, for all $z\in L^{\infty}X\rtimes \Gamma$, or, equivalently, that $\lim_{n\rightarrow\infty}||\phi_n(f)-f||_2=0,$ for every $f\in L^{\infty}X$. By approximating $f$ (in $||.||_2$) with continuous functions, it suffices to show that if $g\in C(X)$, then $\lim_{n\rightarrow\infty}||\phi_n(g)-g||_{\infty}=0$. This in turn follows by using the fact that $g$ is uniformly continuous together with the following estimate
 $$ ||\phi_n(g)-g||_{\infty}=\sup_{x\in X}|\int_{X}\eta_n(x,y)(g(y)-g(x))d\mu(y)|\leq \sup_{d(x,y)<1/n}|g(y)-g(x)|.$$ 

Finally, 
to get the conclusion, we only need to show that the $\Phi_n$'s are compact relative to $L\Gamma$. 
For this, fix $n\geq 1$ and $\varepsilon>0$. Next, denote $N=L^{\infty}X\rtimes \Gamma$ and identify the Hilbert space $L^2N$ with $L^2X\overline{\otimes}l^2\Gamma$, in the natural way. Then, in this identification, we have that $(<N,e_{L\Gamma}>,Tr)=(\Bbb B(L^2X)\overline{\otimes}L\Gamma,tr\otimes\tau),$ where $tr$ is the natural trace on $\Bbb B(L^2X)$,  and  $T_{\Phi_n}=T_{\phi_n}\otimes 1$. Since $T_{\phi_n}\in\Bbb B(L^2X)$ is a compact operator (being  Hilbert-Schmidt by definition) we can find a finite dimensional projection  $p\in\Bbb B(L^2X)$ such that 
$||T_{\phi_n}(1-p)||\leq\varepsilon$. Thus, if we let $P=p\otimes 1$, then $P$ has finite trace and $||T_{\Phi_n}(1-P)||= ||T_{\phi_n}(1-p)||\leq\varepsilon$.
\vskip 0.05in
$(b)\Longrightarrow (c)$ This is a consequence of Proposition 3.4. in [P6].
\vskip 0.05in
$(c)\Longrightarrow (a)$ This implication follows from Proposition 6.10. below.\hfill$\blacksquare$
\vskip 0.1in
Before proceeding to the last result of this Section, let us observe an immediate corollary of Theorem 6.9.:
 if $\Gamma\curvearrowright^{\sigma}(X,\mu)$ is a compact action of a group $\Gamma$ with Haagerup's property then $L^{\infty}X\rtimes_{\sigma}\Gamma$ has Haggerup's property.  Assuming that $\Gamma$ has Haagerup's property, let $\psi_n:\Gamma\rightarrow\Bbb C$ be a sequence of positive definite functions such that $\lim_{n\rightarrow\infty}\psi_n(\gamma)=1$, for all $\gamma\in\Gamma$, $\psi_n(e)=1$ and $\lim_{\gamma\rightarrow\infty}\psi_n(\gamma)=0$, for all $n$. Then, in the notation of the above proof, it is easy to check that $\chi_n(\sum_{\gamma}f_{\gamma}u_{\gamma})=\sum_{\gamma}\psi_n(\gamma)\phi_n(f_{\gamma})u_{\gamma}$ defines a sequence of unital completely positive compact maps on $L^{\infty}X\rtimes_{\sigma}\Gamma$ which satisfy $\lim_{n\rightarrow}||\chi_n(z)-z||_2=0$, for all $z\in L^{\infty}X\rtimes_{\sigma}\Gamma$.  In particular, this observation gives a different proof of Corollary 3.5. in [Jo2].

Note that the next proposition has been obtained in [NS] under the assumption that $\Gamma$ is abelian. 
\vskip 0.1in 
\proclaim {6.10. Proposition} Let $\Gamma\curvearrowright^{\sigma}(X,\mu)$ be an ergodic measure preserving action. Denote by $M=L^{\infty}X\rtimes_{\sigma}\Gamma$ the crossed product von Neumann algebra associated to $\sigma$ and by $L\Gamma$ the von Neumann subalgebra generated by the canonical unitaries $\{u_{\gamma}|\gamma\in\Gamma\}$. Let $\Gamma\curvearrowright^{\alpha}(Y,\nu)$ be the maximal compact quotient  of $\sigma$. 

Then $q\Cal N_{M}(L\Gamma)''=L^{\infty}Y\rtimes_{\alpha}\Gamma.$
In particular, $\sigma$ is compact iff $L\Gamma$ is quasi-regular in $M$ and $\sigma$ is weakly mixing iff $q\Cal N_{M}(L\Gamma)=L\Gamma$.
\endproclaim
{\it Proof.}  Let $a\in L^{\infty}Y$ such that the linear span $\Cal H$ of $\{\gamma a|\gamma\in\Gamma\}$ is finite dimensional. Since $aL\Gamma\subset L\Gamma\Cal H$ and $L\Gamma a\subset\Cal H L\Gamma$, we get that $a\in q\Cal N_{M}(L\Gamma)$.
As the set of such  $a$'s in $||.||_2$-dense in $L^{\infty}Y$ ($\alpha$ is compact), we deduce that $L^{\infty}Y\subset q\Cal N_{M}(L\Gamma)''.$  
Since we also have that $L\Gamma \subset q\Cal N_{M}(L\Gamma)$, altogether we get that $L^{\infty}Y\rtimes_{\alpha}\Gamma\subset q\Cal N_{M}(L\Gamma)''$.

Now, if $x\in q\Cal N_{M}(L\Gamma)$, then $E_{L^{\infty}Y\rtimes_{\alpha}\Gamma}(x)\in q\Cal N_{M}(L\Gamma)$.
 Thus,  in order to derive the reverse inclusion it suffices to show that if $x\in M$ satisfies $x\perp  L^{\infty}Y\rtimes_{\alpha}\Gamma$ and $L\Gamma x\subset\sum_{i=1}^n x_iL\Gamma$, for some  $x_1,x_2,..,x_n\in M$, then $x=0$. Note that the proof of this fact is based on an argument due to S. Popa (see section 3 in [P7] and the proofs of 1.1. in [IPP] and 2.10. in [P7]).
\vskip 0.02in
For every $\varepsilon>0$, let $f_{\varepsilon}(t)=1_{(\varepsilon,+\infty)}(t)t^{-1/2}$ and define $x_{\varepsilon}=xf_{\varepsilon}(E_{L\Gamma}(x^*x))$.
Then for all $\varepsilon>0$ we have that  $L\Gamma x_{\varepsilon}\subset \sum x_iL\Gamma$, $x_{\varepsilon}\perp L^{\infty}Y\rtimes_{\alpha}\Gamma$ and $E_{L\Gamma}(x_{\varepsilon}^*x_{\varepsilon})$ is a projection. Moreover, if $x_{\varepsilon}=0$, for all $\varepsilon>0$, then $x=0$. 
From all this we deduce that we can assume that $q:=E_{L\Gamma}(x^*x)$ is a projection, or, equivalently, that the orthogonal projection onto the closure of $xL\Gamma$ equals $xe_{L\Gamma}x^*$. 
Denote by $p$ the orthogonal projection onto the closure of $L\Gamma x L\Gamma$. Then we have that $p\in <M,e_{L\Gamma}>$, Tr$(p)<\infty$, $p\leq 1-e_{L\Gamma}$ and  $q_{\gamma}:=u_{\gamma}xe_{L\Gamma}x^*u_{\gamma}^*\leq p,$ for all $\gamma\in\Gamma.$

Now, let $\{1=\eta_0,\eta_1,..,\eta_n,..\}\subset M$ be an orthonormal basis of $M$ over $L\Gamma$ and for all $n\geq 1$, denote $f_n=\sum_{i=1}^n \eta_ie_{L\Gamma}{\eta_i}^*$.
Using  the Cauchy-Schwarz inequality we get that for all $\gamma\in\Gamma$ and every $n$, $$\text{Tr}(q)=\text{Tr}(q_{\gamma})=\text{Tr}(q_{\gamma}p)\leq\tag 6.4.a.$$ $$\text{Tr}(f_nq_{\gamma}f_np)+2||q_{\gamma}||_{2,\text{Tr}}||f_np-p||_{2,\text{Tr}}\leq\text{Tr}(f_nq_{\gamma})+2||q||_{2,\text{Tr}}||f_np-p||_{2,\text{Tr}},$$
where by definition $||x||_{2,Tr}=Tr(x^*x)^{1/2}$.

We next claim that for every $\varepsilon>0$ and for all $n\geq 1$, we can find $\gamma\in\Gamma$ such that $\text{Tr}(f_nq_{\gamma})\leq\varepsilon.$
 Note first that $$\text{Tr}(f_nq_{\gamma})=\sum_{i=1}^n\text{Tr}(\eta_ie_{L\Gamma}{\eta_i}^*u_{\gamma}xe_{L\Gamma}x^*{u_{\gamma}}^*)=\sum_{i=1}^n||E_{L\Gamma}({\eta_i}^*u_{\gamma}x)||_2^2\tag 6.4.b.$$

Since $x,\eta_1,..,\eta_n\perp L^{\infty}Y\rtimes_{\alpha}\Gamma$, by approximating $x$ and the $\eta_i$'s with finitely supported vectors (i.e. vectors of the form $\sum a_{\gamma}u_{\gamma}$ with $a_{\gamma}=0$, outside a finite set $F\subset\Gamma$) we can find $a_1,..,a_k\in L^2X\ominus L^2Y$ such that $$\sum_{i=1}^n||E_{L\Gamma}({\eta_i}^*u_{\gamma}x)||_2^2\leq\varepsilon/2+\sum_{j,l=1}^k|\tau(a_i\gamma(a_j))|^2,\forall\gamma\in\Gamma\tag 6.4.c.$$ Since $\alpha$ is the maximal compact quotient of $\sigma$ and since $a_i\perp L^2Y$, for all $i$, we can find $\gamma\in\Gamma$ such that $\sum_{j,l=1}^k|\tau(a_i\gamma(a_j))|^2\leq\varepsilon/2$. By combining this inequality with (6.4.b.) and (6.4.c.), we get the claim. 

To finish the proof of the proposition, let $\varepsilon>0$. Since $p$ has finite trace and $p\leq 1-e_{L\Gamma}$, we can find $n$ such that $||f_np-p||_{2,\text{Tr}}\leq \varepsilon (4||q||_{2,\text{Tr}}+1)^{-1}$. The above claim gives an element $\gamma\in\Gamma$ such that $\text{Tr}(f_nq_{\gamma})\leq\varepsilon/2.$ Finally, by using these inequalities together with (6.4.b.), we get that $\text{Tr}(q)\leq\varepsilon$. Since $\varepsilon>0$ is arbitrary, it follows that $q=0$, thus $x=0$.\hfill$\blacksquare$

\vskip 0.1in
{\it Acknowledgment.} I am grateful to Professor Sorin Popa for helpful suggestions.
\vskip 0.1in

\head References\endhead
\item {[Bu]} M. Burger: {\it Kazhdan constants for SL$(3,\Bbb Z)$}, J. Reine Angew. Math. {\bf 413}(1991), 36--67.
\item {[CFW]} A. Connes, J. Feldman, B. Weiss: {\it  An amenable equivalence relation is generated by a single transformation},  Ergodic Theory Dynamical Systems {\bf 1}(1981), no. 4, 431--450.
\item {[Co]} S. Coskey: {\it Descriptive aspects of torsion-free abelian groups}, PhD Thesis, Rutgers University, 2008.
\item {[Dy]} H. Dye: {\it On groups of measure preserving transformations,} I, Amer. Math. J. {\bf 81}(1959), 119--159.
\item {[ET]} I. Epstein, T. Tsankov: {\it  Modular actions and amenable representations,} preprint 2007, to appear in Trans. Am. Math. Soc.
\item {[FM]} J.Feldman, C.C. Moore: {\it Ergodic equivalence relations, cohomology, and von Neumann algebras, II}, Trans. Am. Math. Soc. {\bf 234}(1977), 325--359.
\item {[Fu1]} A. Furman, {\it Gromov's measure equivalence and rigidity of higher rank lattices},  Ann. of Math. (2) {\bf 150}(1999),  1059--1081.
\item {[Fu2]} A. Furman: {\it Orbit equivalence rigidity,} Ann. of Math. (2)  {\bf 150}(1999), 1083--1108.
\item {[Fu3]} A. Furman: {\it Outer automorphism groups of some ergodic equivalence relations,}  Comment. Math. Helv.  {\bf 80} (2005), 157--196.
\item {[Fu4]} A. Furman: {\it On Popa's Cocycle Superrigidity Theorem}, IMRN  {\bf 2007},  no. 19, Art. ID rnm073, 46 pp.
\item {[G]} E. Glasner: {\it Ergodic theory via joinings}, Mathematical Surveys and Monographs, 101. AMS, Providence, RI, 2003. xii+384 pp.
\item{[Hj]} G. Hjorth: {\it A converse to Dye's theorem}, Trans. Am. Math. Soc. {\bf 357}(2005), 3083-3103.
\item {[IPP]} A. Ioana, J. Peterson, S. Popa: {\it Amalgamated free products of weakly rigid factors
and calculation of their symmetry groups}, Acta Math. {\bf 200}(2008), no. 1, 85--153.
\item {[Jo1]} P. Jolissaint: {\it On the property (T) for pairs of topological groups,} Enseign. Math.  {\bf 51}(2005), no. 1-2, 31--45.
\item {[Jo2]} P. Jolissaint: {\it Actions of dense subgroups of compact groups and II$_1$-factors with the Haagerup property,}  Ergodic Theory Dynam. Systems  {\bf 27}(2007), no. 3, 813--826.
\item {[K]} D. Kazhdan: {\it On the connection of the dual space of a group with the structure of its closed subgroups}, Funct. Anal. and its Appl. {\bf 1}(1967), 63--65.
\item {[Ki1]} Y. Kida: {\it Measure equivalence rigidity of the mapping class group,} 
preprint math math.GR/0607600, to appear in
Ann. Math.
\item {[Ki2]} Y. Kida: {\it  Orbit equivalence rigidity for ergodic actions of the mapping
 class group},  Geom. Dedicata  {\bf 131} (2008), 99--109.
\item {[Ma]} G. Margulis: {\it Finitely-additive invariant measures on Euclidian spaces,}  Ergodic Theory Dynam. Systems {\bf 2}(1982), 383--396.
\item {[MSh1]} N. Monod, Y. Shalom: {\it  Orbit equivalence rigidity and bounded cohomology}, Ann. Math.
 {\bf 164} (2006) 825–-878.
\item {[MSh2]} N. Monod, Y. Shalom: {\it Cocycle superrigidity and bounded cohomology for negatively curved spaces,}  J. Differential Geom.  {\bf 67} (2004),  no. 3, 395--455.
\item {[MvN]} F. Murray , J. von Neumann: {\it On rings of operators,} Ann. Math. {\bf 37} (1936), 116–229.
\item {[NS]} S. Neshveyev, E. St{\o}rmer: {\it Ergodic Theory and Maximal Abelian Subalgebras of the Hyperfinite Factor},  J. Funct. Anal. {\bf 195} (2002), 239--261.
\item {[OP]} N. Ozawa, S. Popa: {\it On a class of II$_1$ factors with at most one Cartan subalgebra}, math.OA/0706.3623, to appear in Ann. Math.
\item {[OW]} D. Ornstein, B.Weiss: {\it Ergodic theory of amenable groups. I. The Rokhlin lemma.}, Bull. Amer. Math. Soc. (N.S.) {\bf 1} (1980), 161-164.
\item {[P1]} S. Popa: {\it Strong Rigidity of II$_1$ Factors Arising from Malleable Actions of w-Rigid Groups II}, Invent. Math. {\bf 165} (2006), 409--451.
\item {[P2]} S. Popa: {\it Cocycle and orbit equivalence superrigidity for malleable actions 
of $w$-rigid groups,}  Invent. Math.  {\bf 170}(2007),  no. 2, 243--295.
\item {[P3]} S. Popa: {\it On the superrigidity of malleable actions with spectral gap},  
math.GR/0608429, to appear in Journal of the AMS.
\item {[P4]} S. Popa: {Deformation and rigidity for group actions and von Neumann algebras},  International
 Congress of Mathematicians. Vol. I,  445--477, Eur. Math. Soc., Z$\ddot{u}$rich, 2007. 
\item {[P5]} S. Popa: {\it Some rigidity results for non-commutative Bernoulli shifts,}  J. Funct. Anal.  {\bf 230}  (2006),  273--328.
\item {[P6]} S. Popa: {\it On a class of type II$_1$ factors with Betti numbers invariants}, Ann. Math. {\bf 163}(2006), 809--889. 
\item {[P7]} S. Popa: {\it Strong rigidity of II$_1$ factors arising from malleable actions of $w$-rigid groups I}, 
  Invent. Math.  {\bf 165} (2006),  369--408.
\item {[PV]} S. Popa, S. Vaes: {\it  Strong rigidity of generalized Bernoulli actions and computations of their 
symmetry groups.} Adv. Math. {\bf 217}(2008),  no. 2, 833--872.
\item {[Sh1]} Y. Shalom: {\it Measurable group theory,}  European Congress of Mathematics,  391--423, Eur. Math. Soc., Z$\ddot u$rich, 2005. 
\item {[Sh2]} Y. Shalom: {\it Rigidity of commensurators and irreducible lattices},  Invent. Math.  {\bf 141}(2000),  1--54.
\item {[T]} S. Thomas: {\it Superrigidity and countable Borel equivalence relations},  Ann. Pure Appl. Logic  {\bf 120}(2003),  no. 1-3, 237--262.
\item {[Z]} R. Zimmer: {\it Ergodic theory and semisimple groups}
Monographs in Mathematics, 81. Birkhäuser Verlag, Basel, 1984. x+209 pp.
\enddocument